\documentclass{amsart}
\usepackage{amssymb,amscd,amsthm,verbatim}
\input epsf                 
\allowdisplaybreaks

\newtheorem{proposition}{Proposition}[section]
\newtheorem{theorem}[proposition]{Theorem}
\newtheorem{corollary}[proposition]{Corollary}
\newtheorem{lemma}[proposition]{Lemma}

\newtheorem{prop}[proposition]{Proposition}

\newcommand{\reals}{\mathbb R}

\newcommand{\supp}{\mathrm{supp}}

\newcommand{\Sh}{{\mathcal Sh}}
\newcommand{\Q}{{\mathcal Q}}
\newcommand{\A}{{\mathcal A}}
\newcommand{\R}{{\mathcal R}}

\newcommand{\B}{{\mathcal {B}}}
\newcommand{\F}{{\mathcal F}}
\newcommand{\Po}{{\mathcal P}}

\newcommand{\covers}{{\,\,\,\cdot\!\!\!\! >\,\,}}
\newcommand{\covered}{{\,\,<\!\!\!\!\cdot\,\,\,}}

\newcommand{\set}[1]{{\lbrace #1 \rbrace}}

\newcommand{\join}{\vee}
\newcommand{\bigjoin}{\bigvee}
\newcommand{\meet}{\wedge}

\newcommand{\Cg}{\mathrm{Cg}}
\newcommand{\Con}{\mathrm{Con}}
\newcommand{\relint}{\mathrm{relint}}
\newcommand{\Span}{\mathrm{Span}}
\newcommand{\Irr}{\mathrm{Irr}}

\newcommand{\pidown}{\pi_\downarrow}
\newcommand{\piup}{\pi^\uparrow}

\newcommand{\case}[2]{\item[Case #1:] #2 \\ }

\newcommand{\br}[1]{\langle #1 \rangle}
\newcommand{\ul}[1]{\underline{#1}}

\begin{document}
\title[Lattice congruences]{Lattice congruences of the weak order}

\author{Nathan Reading}
\address{
Mathematics Department\\
       University of Michigan\\
       Ann Arbor, MI 48109-1109\\
USA}
\thanks{The author was partially supported by NSF grant DMS-0202430.}
\email{nreading@umich.edu}
\urladdr{http://www.math.lsa.umich.edu/$\sim$nreading/}

\begin{abstract}
We study the congruence lattice of the poset of regions of a hyperplane arrangement, with particular emphasis on the weak order on a finite
Coxeter group.
Our starting point is a theorem from a previous paper which gives a geometric description of the poset of join-irreducibles of the congruence 
lattice of the poset of regions in terms of certain polyhedral decompositions of the hyperplanes.
For a finite Coxeter system $(W,S)$ and a subset $K\subseteq S$, let $\eta_K:w\mapsto w_K$  be the projection onto the parabolic subgroup $W_K$.
We show that the fibers of $\eta_K$ constitute the smallest lattice congruence with $1\equiv s$ for every $s\in(S-K)$.
We give an algorithm for determining the congruence lattice of the weak order for any finite Coxeter group and for a finite Coxeter group of 
type~A or~B we define a directed graph on subsets or signed subsets such that the transitive closure of the directed graph is the poset of 
join-irreducibles of the congruence lattice of the weak order.
\end{abstract}

\maketitle

\section{Introduction}
\label{intro}

A {\em congruence} in a lattice $L$ is an equivalence relation on the elements of $L$ which respects the join and meet operations in $L$.
A {\em join-irreducible} in $L$ is an element $\gamma$ which covers a unique element $\gamma_*$.
When $L$ is a finite lattice, congruences are determined by the set of join-irreducibles such that $\gamma\equiv \gamma_*$.
The set of congruences of a finite lattice, partially ordered by refinement, forms a finite distributive lattice~\cite{Fun-Nak} denoted $\Con(L)$.
In general it may not be straightforward to determine $\Con(L)$ for a given lattice $L$, and indeed there are few examples
where the congruence lattices are known for a class of finite lattices.
Exceptions include the finite distributive lattices, whose congruence lattice is a Boolean algebra and the Tamari lattice $T_n$, whose 
congruence lattice was determined by Geyer~\cite{Geyer}.
One construction of the Tamari lattice is as the subposet (in fact, sublattice~\cite{Nonpure II}) of the weak order on permutations consisting of 
$312$-avoiding permutations.

A finite lattice is called {\em congruence uniform} (or sometimes {\em bounded}) when the join-irreducibles of $\Con(L)$ are in 
bijection with the join-irreducibles of $L$ and when the dual condition holds.
In this case, the subposet $\Irr(\Con(L))$ consisting of join-irreducibles can be thought of as a partial order on the join-irreducibles of $L$.
Geyer showed that the Tamari lattice is congruence uniform and determined $\Irr(\Con(T_n))$ as a partial order on the join-irreducibles
of $T_n$.
Caspard, Le Conte de Poly-Barbut and Morvan~\cite{boundedref} showed that the weak order on a finite Coxeter group is congruence 
uniform.
This result was generalized to the posets of regions of certain hyperplane arrangements as \cite[Theorem 25]{hyperplane}
which also provided an explicit description of the congruence lattice.
This description involved cutting the hyperplanes into pieces called {\em shards} and using incidence relations among the shards to define a 
directed graph whose transitive closure is the poset of irreducibles of the congruence lattice.

Using this geometric description, we give an algorithm for determining the congruence lattice of the weak order on a finite Coxeter group 
and prove, in a more general context, the result on parabolic subgroups described in the abstract, as well as some related results.
We give a concrete description of the congruence lattices of the weak orders on finite Coxeter groups of types A and B\@.
Specifically, we reduce the incidence relations among the shards to combinatorial criteria on subsets and signed subsets.
The combinatorics of the congruence lattice in type A is very simple, but in type B it is significantly more complicated.
The concrete description of the congruence lattice has important applications.
It is used in~\cite{Cambrian} to make a broad generalization of the Tamari lattices, called the {\em Cambrian lattices}.
For each finite Coxeter group the Cambrian lattices are a family of lattices, defined as quotients of the weak order with respect to certain
natural congruences.
The Tamari lattice $T_n$ is one of the Cambrian lattices associated to the Coxeter group $S_n$.
The congruence lattices of the Cambrian lattices are also determined, generalizing a result of Geyer~\cite{Geyer} on the Tamari lattice. 
The concrete description is also used in~\cite{con_app} to show that certain congruences of the weak order on $S_n$ give rise to sub Hopf 
algebras of the Malvenuto-Reutenauer Hopf algebra of permutations~\cite{MR}.
Also in~\cite{con_app}, for any lattice quotient of the weak order, a fan is constructed such that the geometry of the fan interacts with the 
lattice quotient in many of the pleasant ways that the weak order interacts with the geometry of the associated hyperplane arrangement.
Finally, the algorithm and the explicit description lend themselves to using a computer to count the numbers of congruences of the weak orders.

This paper is the first in a series of papers which continues in~\cite{con_app} and~\cite{Cambrian}.
Each paper relies on the results of the preceding papers and cites later papers only for motivation or in the context of examples.

This paper is organized as follows:
Section~\ref{lattice} introduces background and notation concerning lattices and the congruence lattice.
In Section~\ref{poset} we provide necessary background on hyperplane arrangements, the poset of regions and shards.
In Section~\ref{cox} we give background information on Coxeter groups and describe an algorithm for determining the congruence lattice of 
the weak order, and in Section~\ref{coxab} we give standard combinatorial interpretations of the Coxeter groups of types~A and~B\@.
Section~\ref{first} deals with parabolic congruences of the poset of regions in the case when it is a lattice.
Sections~\ref{Con Bn} and~\ref{Con Sn} determine the shard digraph $\Sh(W)$, whose transitive closure is $\Irr(\Con(W))$, where 
$W$ refers to the weak order on a Coxeter group of type B or A (the types are treated in that order).
We conclude with some examples and applications in Section~\ref{examples}.

\section{The congruence lattice of a lattice}
\label{lattice}

Our notation for posets and lattices is fairly standard.
Some pieces of notation that may not be familiar are as follows:
We use the symbol ``\!\!$\covered$\!\!'' to denote cover relations in a poset.
If $A$ and $B$ are sets, we use the notation $A-B$ to mean the set of elements in $A$ but not in $B$.
This notation is not meant to imply that $B\subseteq A$.

Let $L$ be a finite lattice, with join and meet operations denoted by $\join$ and $\meet$ respectively.
An element $\gamma$ of $L$ is join-irreducible if for any $S\subseteq L$ with $\gamma=\join S$, we have $\gamma\in S$.
Equivalently, $\gamma$ is join-irreducible if it covers exactly one element $\gamma_*$ of $L$.
Meet-irreducible elements are defined dually, and a meet-irreducible $\mu$ is covered by a unique element $\mu^*$.
Denote the set of join-irreducibles of $L$ by $\Irr(L)$.
We also use $\Irr(L)$ to represent the induced subposet of $L$ consisting of join-irreducible elements.  
A {\em congruence} of $L$ is an equivalence relation on the elements of $L$ which respects joins and meets in the following sense:
If $a_1\equiv a_2$ and $b_1\equiv b_2$ then $a_1\join b_1\equiv a_2\join b_2$ and similarly for meets.
If $\Theta$ is a congruence such that $x\equiv y$ for some $x\covered y$, we say $\Theta$ {\em contracts} the edge $x\covered y$.
For an element $y$, if there exists an edge $x\covered y$ contracted by $\Theta$, we say $\Theta$ contracts $y$.
In particular $\Theta$ contracts a join-irreducible $\gamma$ if and only if $\gamma\equiv \gamma_*$.
A lattice congruence is determined by the set of join-irreducibles it contracts (see for example \cite[Section II.3]{Free Lattices}).

Lattice congruences can be described entirely in terms of the partial order on $L$, without referring to joins and meets.
This leads to the definition of an {\em order congruence} on a poset $P$.
A general definition can be found in~\cite{Cha-Sn}, but here we restrict to the finite case, and we do not follow the convention 
of~\cite{Cha-Sn} that the equivalence relation with only one class is always a congruence.

Let $P$ be a finite poset with an equivalence relation $\Theta$ defined on the elements of $P$.
Given $a\in P$, let $[a]_\Theta$ denote the equivalence class of $a$.
The equivalence relation is an {\em order congruence} if:
\begin{enumerate}
\item[(i) ] Every equivalence class is an interval.
\item[(ii) ] The projection $\pidown:P\rightarrow P$, mapping each element $a$ of $P$ to the minimal element in $[a]_\Theta$, is order-preserving.
\item[(iii) ] The projection $\piup:P\rightarrow P$, mapping each element $a$ of $P$ to the maximal element in $[a]_\Theta$, is 
order-preserving.
\end{enumerate}
If $P$ is a lattice, then $\Theta$ is an order congruence if and only if it is a lattice congruence.

Define a partial order on the congruence classes by $[a]_\Theta\le[b]_\Theta$ if and only if there exists $x\in[a]_\Theta$ and $y\in[b]_\Theta$ 
such that $x\le_Py$.
The set of equivalence classes under this partial order is $P/\Theta$, the {\em quotient} of $P$ with respect to $\Theta$.
The quotient $P/\Theta$ is isomorphic to the induced subposet $\pidown(P)$.
The map $\piup$ maps $\pidown(P)$ isomorphically onto $\piup(P)$.
The inverse is $\pidown$.
When $P$ is a lattice, this construction of the quotient corresponds to the algebraic construction of the quotient lattice mod $\Theta$.
For more information on order congruences and quotients, see~\cite{Cha-Sn,Order}.

Congruences on $L$ are, in particular, partitions of the elements of $L$, and $\Con(L)$ is the set of congruences of $L$
partially ordered by refinement.
The partial order $\Con(L)$ is a distributive lattice, and thus is uniquely determined by the subposet $\Irr(\Con(L))$.
The meet in $\Con(L)$ is the intersection of the congruences as relations.
If $\Theta_1$ and $\Theta_2$ are congruences on $L$, with associated downward projections $(\pidown)_1$ and $(\pidown)_2$, let 
$\Theta_1\join\Theta_2$ have associated downward projection $\pidown$.
Then $x\in L$ has $\pidown x=x$ if and only if both $(\pidown)_1 x=x$ and $(\pidown)_2 x =x$.
Thus the quotient of $L$ mod $\Theta_1\join\Theta_2$ is isomorphic to the induced subposet $((\pidown)_1L)\cap((\pidown)_2L)$ of $L$.

Given lattices $L_1$ and $L_2$ a {\em homomorphism} from $L_1$ to $L_2$ is a map $\eta:L_1\to L_2$ such that for all $x$ and $y$ in 
$L_1$ we have $\eta(x\join y)=\eta(x)\join\eta(y)$ and similarly for meets.
Given a lattice homomorphism $\eta$, the equivalence relation whose classes are the fibers of $\eta$ is a congruence,
and conversely, given a congruence $\Theta$ on $L$, the map from an element to its equivalence class is a homomorphism $L\to(L/\Theta)$.
Alternately, the map $\pidown$ is a homomorphism from $L$ to $\pidown(L)\cong L/\Theta$.
If $\eta_1:L\to L_1$ and $\eta_2:L\to L_2$ are lattice homomorphisms, we say $\eta_2$ {\em factors through} $\eta_1$ if there is a lattice
homomorphism $\eta:L_1\to L_2$ such that $\eta_2=\eta\circ\eta_1$.
Let $\Theta_1$ and $\Theta_2$ be the lattice congruences associated to $\eta_1$ and $\eta_2$.
If $\Theta_1\le\Theta_2$ in $\Con(L)$ then $\eta_2$ factors through $\eta_1$.

Given a covering pair $x\covered y$ in $L$, let $\Cg(x,y)$ be the smallest lattice congruence contracting that edge.
Then $\Cg(x,y)$ is a join-irreducible congruence.
Given a join-irreducible $\gamma$ of $L$, write $\Cg(\gamma)$ for $\Cg(\gamma,\gamma_*)$.
The map $\Cg:\Irr(L)\to\Irr(\Con(L))$ is onto, but need not be one-to-one.
A lattice $L$ is {\em join congruence uniform} if $\Cg$ is a bijection~\cite{cong norm}.
The notion of {\em meet congruence uniformity} is defined dually, and $L$ is called {\em congruence uniform} if it is both join and meet 
congruence uniform.
All of the lattices we deal with here are self-dual, so that congruence uniformity and join congruence uniformity coincide.
When $L$ is a congruence uniform lattice, $\Irr(\Con(L))$ can be thought of as a partial order on the join-irreducibles.
If $\Theta$ is a congruence on a congruence uniform lattice $L$, then $\Irr(\Con(L/\Theta))$ is the filter in $\Irr(\Con(L))$ consisting 
of join-irreducibles of $L$ not contracted by $\Theta$.

A lattice $L$ is called {\em join semi-distributive} if for any $x,y,z\in L$
\[x\join y=x\join z \mbox{ implies }x\join(y\meet z)=x\join y.\]
A lattice is {\em meet semi-distributive} if the dual condition holds and {\em semi-distributive} if it is both join semi-distributive and meet 
semi-distributive.
Congruence uniform lattices are in particular {\em semi-distributive}~\cite{cong norm}.

A finite lattice is join-semi-distributive if and only if for every join-irreducible~$\gamma$ there is a unique meet-irreducible~$\mu$ such that 
$\gamma\meet\mu=\gamma_*$ and $\gamma\join\mu=\mu^*$~\cite[Theorem~2.56]{Free Lattices}, and meet-semi-distributive if and only if the dual condition holds.
In~\cite{hyperplane}, such a pair $(\gamma,\mu)$ is called a {\em subcritical pair}. 
There it is also shown~\cite[Proposition~17]{hyperplane} that if $\Theta$ is a congruence on $L$ and $(\gamma,\mu)$ is a subcritical pair in 
$\pidown(L)$, then $(\gamma,\piup(\mu))$ is a subcritical pair in $L$.

If $(\gamma,\mu)$ is a subcritical pair, then a congruence $\Theta$ contracts $\gamma$ if and only if it contracts the edge $\mu\covered\mu^*$.
Let $L$ be semi-distributive and let $\sigma$ be the map taking each $\gamma$ to the unique $\mu$ such that $(\gamma,\mu)$ is a subcritical pair.
If $\alpha$ is an anti-automorphism of $L$, then the composition $\alpha\circ\sigma$ is an automorphism of $\Irr(\Con(L))$.
Suppose $\Theta_1$ and $\Theta_2$ are two congruences of $L$, contracting respectively the sets $C_1$ and $C_2:=\alpha(\sigma(C_1))$ 
of join-irreducibles.
Then $\sigma$ maps $C_1$ to  the set of meet-irreducibles $\mu$ such that $\mu\covered\mu^*$ is contracted by $\Theta_1$.
The fact that $C_2$ is the image of that set under $\alpha$ means that $\alpha$ maps $\Theta_1$-classes to $\Theta_2$-classes.
Thus $\alpha$ induces an anti-isomorphism from $L/\Theta_1$ to $L/\Theta_2$.
For more information on $\Con(L)$, congruence uniformity and semi-distributivity, see \cite{Free Lattices,Gratzer}.

\section{The poset of regions}
\label{poset}
In this section, we provide the necessary background on hyperplane arrangements and their posets of regions.
We define the shards of an arrangement and the shard digraph, and determine some of the properties of the shard digraph.

A {\em hyperplane arrangement} $\A$ is a finite nonempty collection of {\em hyperplanes} (codimension~1 subspaces) in $\reals^d$.
In general, the hyperplanes can be affine subspaces, but in this paper we assume all arrangements to be central, meaning that all hyperplanes
contain the origin.
The complement of the union of the hyperplanes is disconnected, and the closures of its connected components are called {\em regions}.
The {\em span} of~$\A$, written $\Span(\A)$, is understood to mean the linear span of the normal vectors of~$\A$, and the {\em rank} 
of~$\A$ is the dimension of $\Span(\A)$.
A region is called {\em simplicial} if the collection of normal vectors to its facets is linearly independent, and a central hyperplane arrangement 
is called {\em simplicial} if every region is simplicial.

The {\em poset $\Po(\A,B)$ of regions} of $\A$ with respect to a fixed region $B$ is a partial order whose elements are the regions, with order 
relations given as follows.
For convenience we think of each hyperplane as a non-specified linear functional which is zero on $H$ and negative in $B$.
For any region $R$, and $x$ in the interior of $R$, the set $S(R):=\set{H\in\A:H(x)>0}$ is independent of the choice of $x$ and is called
the {\em separating set} of $R$.
The poset of regions is a partial order on the regions with $R_1\le R_2$ if and only if $S(R_1)\subseteq S(R_2)$.
The fixed region $B$, called the {\em base region}, is the unique minimal element.
The antipodal map induces an anti-automorphism on $\Po(\A,B)$, called the {\em antipodal anti-automorphism},  denoted by $R\mapsto -R$.
This is an anti-automorphism because $S(-R)=\A-S(R)$.
In particular, there is a unique maximal element $-B$, which is the minimal upper bound of the set of atoms (elements covering $B$) of
$\Po(\A,B)$.
For more details on hyperplane arrangements and the poset of regions, see~\cite{BEZ,Edelman,OT,hyperplane,hplanedim}.

Let $\A$ be an arrangement in $\reals^d$ with base region $B$.
We say $\A'$ is a {\em rank-two subarrangement} if $\A'$ consists of all the hyperplanes containing some subspace of dimension $d-2$ 
and $|\A'|\ge 2$.
There is a unique region $B'$ of $\A'$ containing $B$.
The two facet hyperplanes of $B'$ are called {\em basic} hyperplanes in $\A'$.
If $H_1$ and $H_2$ are basic in $\A'$ but $H\in\A'$ is not, then $(H\cap B') = (H_1\cap H_2\cap B')$.
Intersecting both sides of the equality with~$B$, we obtain the following, which we name as a lemma for easy reference later.
\begin{lemma}
\label{basic containment}
If $H_1$ and $H_2$ are basic in $\A'$ but $H\in\A'$ is not, then $(H\cap B) = (H_1\cap H_2\cap B)$.
\end{lemma}

We define a directed graph $\Q(\A,B)$ whose vertex set is~$\A$, with directed edges $H_1\rightarrow H_2$ whenever $H_1$ 
is basic and $H_2$ is not basic in the rank-two subarrangement determined by $H_1\cap H_2$.
For $H\in\A$, define $a(H):=\set{H':H'\to H\mbox{ in }\Q(\A,B)}$.
It is immediate that $\Q(\A,B)$ contains no two-element cycles.

We cut the hyperplanes of $\A$ into pieces called {\em shards} as follows.
For each non-basic $H$ in a rank-two subarrangement $\A'$, cut $H$ into connected components by removing the subspace $\cap\A'$ from $H$.
Do this cutting for each rank-two subarrangement, and call the closures of the resulting connected components of the hyperplanes {\em shards}.
Note that in \cite{hyperplane,hplanedim}, closures were not taken when the shards were defined.
The difference is entirely one of convenience and does not affect the combinatorics.
In particular, taking closures makes the definition of $\Sh(\A,B)$, below, look different than in the previous papers, but it is the same 
directed graph.
Figures~\ref{Bpic} and~\ref{Apic}, in Sections~\ref{Con Bn} and~\ref{Con Sn}, depict the collection of shards for two different arrangements.

For each shard $\Sigma$, let $H_\Sigma$ be the hyperplane containing~$\Sigma$.
Define the {\em shard digraph} $\Sh(\A,B)$ to be the directed graph whose vertices are the shards, and whose arrows are as follows:
There is a directed arrow $\Sigma_1\rightarrow \Sigma_2$ whenever $H_{\Sigma_1}\to H_{\Sigma_2}$ in $\Q(\A,B)$ and
$\Sigma_1\cap\Sigma_2$ has dimension $d-2$.

Let $U(\Sigma)$ be the set of {\em upper regions} of $\Sigma$, that is, the set of regions $R$ of~$\A$ whose intersection with $\Sigma$ 
has dimension $d-1$ and which have $H_\Sigma\in S(R)$.
Think of $U(\Sigma)$ as an induced subposet of $\Po(\A,B)$.
In this paper regions are denoted by capital letters, so we denote join-irreducible regions by $J$ to emphasize that they are regions.
The following is~\cite[Proposition 2.2]{hplanedim}.
\begin{prop}
\label{j sigma}
A region $J$ is join-irreducible in $\Po(\A,B)$ if and only if $J$ is minimal in $U(\Sigma)$ for some shard $\Sigma$.
\end{prop}

Define $L(\Sigma)$ to be the set of regions $R$ of~$\A$ whose intersection with $\Sigma$ has dimension $d-1$ and which have 
$H_\Sigma\not\in S(R)$.
The map taking a region $R$ to the region $R'$ with $S(R')=S(R)\cup\set{H_\Sigma}$ is an isomorphism from $L(\Sigma)$ to $U(\Sigma)$.
The dual version of Proposition~\ref{j sigma} identifies meet-irreducible regions as maximal elements of $L(\Sigma)$.
The following is~\cite[Proposition 4.1]{hplanedim}.
\begin{prop}
\label{sc sigma}
A pair $(J,M)$ in $\Po(\A,B)$ is subcritical if and only if there is a shard $\Sigma$ such that $J$ is minimal in $U(\Sigma)$, $M$ is 
maximal in $L(\Sigma)$ and $J_*\le M$.
\end{prop}

Since a join-irreducible $J$ covers only one element, there is a {\em unique} shard $\Sigma_J$ associated to $J$.
We write $H_J$ for $H_{\Sigma_J}$.
When $\A$ is simplicial, $\Po(\A,B)$ is a lattice \cite{BEZ}, and furthermore, it is semi-distributive~\cite{hyperplane}.
The following theorem is \cite[Theorem 25]{hyperplane}.
We interpret the transitive closure of $\Sh(\A,B)$ as a poset by setting $\Sigma_1>\Sigma_2$ when $\Sigma_1\to\Sigma_2$.
\begin{theorem}
\label{shard}
Given a simplicial arrangement $\A$, the lattice $\Po(\A,B)$ is congruence uniform if and only if $\Sh(\A,B)$ is acyclic, in which case the 
transitive closure of $\Sh(\A,B)$ is isomorphic to $\Irr(\Con(\Po(\A,B)))$.
\end{theorem}

For a shard $\Sigma$, define $\Phi(\Sigma)$ to be the smallest congruence which contracts any edge $P\covered Q$ such that 
$P\cap Q\subseteq\Sigma$.
In the proof of Theorem 3.3 in~\cite{hyperplane}, independent of the hypotheses that $\A$ is simplicial and $\Sh(\A,B)$ is acyclic, it is shown 
that $\Phi(\Sigma)$ is a well-defined join-irreducible congruence, that $\Phi$ is a surjective map from the shards of $(\A,B)$ to 
$\Irr(\Con(\Po(\A,B)))$ and that $\Sigma_1\to\Sigma_2$ implies $\Phi(\Sigma_1)\ge\Phi(\Sigma_2)$.
Furthermore, although it is not explicitly stated in the proof, it is shown that $U(\Sigma)$ is a connected poset.

\begin{prop}
\label{unique}
If $\A$ is a simplicial arrangement then there is a unique join-irreducible $J_\Sigma$ associated to each shard $\Sigma$.
\end{prop}
\begin{proof}
Suppose for the sake of contradiction that there are two distinct minimal elements of $U(\Sigma)$.
Then since $U(\Sigma)$ is connected, we can find two distinct minimal elements $J_1$ and $J_2$ which are both below the same maximal 
element $M'$ of $U(\Sigma)$.
Let $M$ be the region reached from $M'$ by crossing $\Sigma$.
By Proposition~\ref{sc sigma} $(J_1,M)$ and $(J_2,M)$ are both subcritical pairs, but since $\A$ is simplicial, $\Po(\A,B)$ is 
semi-distributive.
This contradiction shows that $U(\Sigma)$ has a unique minimal element, so by Proposition~\ref{j sigma} we are finished.
\end{proof}
Since $J_\Sigma$ is join-irreducible, it covers a unique region $(J_\Sigma)_*$.
The hyperplane separating $J_{\Sigma}$ from $(J_\Sigma)_*$ is $H_\Sigma$.
A cover relation in $\Po(\A,B)$ is contracted by a congruence $\Theta$ if and 
only if the shard crossed in the covering relation corresponds to a join-irreducible that is contracted by $\Theta$.
In other words, contracting join-irreducibles corresponds to ``removing'' shards.
The ``shard arrangement'' of non-removed shards defines a complete fan $\F_\Theta$ which is studied in depth in~\cite{con_app}.

The following lemmas are helpful in later sections where we apply Theorem~\ref{shard}.
Here ``$\relint$'' means relative interior.
\begin{lemma}
\label{arrows}
There is a directed arrow $\Sigma_1\rightarrow \Sigma_2$ in $\Sh(\A,B)$ if and only if the following three conditions hold:
\begin{enumerate}
\item[(i) ]$H_{\Sigma_1}\to H_{\Sigma_2}$ in $\Q(\A,B)$, 
\item[(ii) ] $H_{\Sigma_1}$ is a facet-defining hyperplane of $\Sigma_2$, and
\item[(iii) ] $(\relint\Sigma_1)\cap\Sigma_2\neq\emptyset$.
\end{enumerate}
\end{lemma}
\begin{proof}
We show that when (i) holds, conditions (ii) and (iii) are equivalent to the requirement that $\Sigma_1\cap\Sigma_2$ has dimension $d-2$.
For convenience let $H_i:=H_{\Sigma_i}$ for $i=1,2$.
Suppose that (i) holds and that $\Sigma_1\cap\Sigma_2$ has dimension $d-2$.
Then $H_1\cap\Sigma_2$ also has dimension $d-2$.
Since $H_1\to H_2$, we have $H_1\cap\,(\relint\Sigma_2)=\emptyset$, which implies (ii).
Since $H_1$ is basic in the rank-two subarrangement determined by $H_1\cap H_2$, the hyperplane $H_2$ does not define a facet of 
$\Sigma_1$.
Thus the fact that $\Sigma_1\cap\Sigma_2$ has dimension $d-2$ implies that $\Sigma_1\cap\Sigma_2$ contains a point in
the relative interior of $\Sigma_1$.

Conversely, suppose (i), (ii) and (iii) hold.
As in the previous paragraph, (i) implies that $H_2$ does not define a facet of $\Sigma_1$.
Let $x\in (\relint\Sigma_1)\cap\Sigma_2$ and let $N$ be a (relatively) open neighborhood of $x$ in $\Sigma_1$.
Since $H_1$ defines a facet of $\Sigma_2$, the intersection $N\cap\Sigma_2$ has dimension $d-2$ and therefore so does $\Sigma_1\cap\Sigma_2$.
\end{proof}

\begin{lemma}
\label{polyhedra}
In a simplicial arrangement $\A$ with base region $B$ the shard $\Sigma$ is the set of points $x\in H_\Sigma$ such that $H'(x)>0$ for 
every $H'\in a(H_\Sigma)\cap S(J_\Sigma)$.
\end{lemma}
\begin{proof}
By definition each point of the relative interior of a facet of $\Sigma$ in $H_\Sigma$ is contained in two hyperplanes 
$H_1,H_2\in a(H_\Sigma)$.
Without loss of generality, $H_1$ separates $\Sigma$ from $B$ and $H_2$ does not.
Thus this facet is defined by the requirement that $H_1$ separate $\Sigma$ from $B$, or alternately by the requirement that 
$H_2$ not separate $\Sigma$ from $B$.
But this means $H_1\in S(J_\Sigma)$ and $H_2\not\in S(J_\Sigma)$.
\end{proof}

\begin{lemma}
\label{Cg}
Let $\A$ be a simplicial arrangement and let $R'\covered R$ in $\Po(\A,B)$ with $S(R)-S(R')=\set{H}$.
Then there is a unique join-irreducible $J$ such that $\Cg(J)=\Cg(R',R)$.
This $J$ is the unique join-irreducible with $H_J=H$ and $a(H_J)\cap S(J)\subseteq S(R)$.
It is also the unique join-irreducible with $H_J=H$ and $J\le R$.
\end{lemma}
\begin{proof}
We have $\Cg(J)=\Cg(R',R)$ if and only if $R'\cap R$ is contained in $\Sigma_J$.
By Proposition~\ref{unique}, for any $\Sigma$ there is a unique $J$ with $\Sigma_J=\Sigma$.
By Lemma~\ref{polyhedra} we have $R'\cap R\subseteq\Sigma_J$ if and only if $H=H_J$ and $H'(x)\ge 0$ for every $x\in R'\cap R$ and
every $H'\in a(H_J)\cap S(J)$.
This is equivalent to the requirement that $H_J=H$ and $a(H_J)\cap S(J)\subseteq S(R)$.
If $R'\cap R\subseteq\Sigma_J$, we have $R\in U(\Sigma_J)$, so $J\le R$.
Conversely, if $J_1\le R$, $H_{J_1}=H$, $J_2\le R$ and $H_{J_2}=H$, then in particular $a(H_{J_1})\cap S(J_1)\subseteq S(R)$ and 
$a(H_{J_2})\cap S(J_2)\subseteq S(R)$, so $J_1=J_2$.
\end{proof}

\begin{lemma}
\label{source}
Let $\A$ be a simplicial arrangement, let $B$ be a base region and let $\Sigma$ be a shard.
Then the following are equivalent:
\begin{enumerate}
\item[(i) ]$\Sigma$ is a source in $\Sh(\A,B)$.
\item[(ii) ]$\Sigma$ is a facet hyperplane of $B$.
\item[(iii) ]There is no facet of $\Sigma$ intersecting $J_\Sigma$ in dimension $d-2$.
\end{enumerate}
\end{lemma}
\begin{proof}
We begin by showing that (i) is equivalent to (ii).
Any shard that arises by cutting a hyperplane is arrowed to by a shard in the hyperplane that cut it, and any shard that is a whole hyperplane 
is not arrowed.
Thus $\Sigma$ is a source in $\Sh(\A,B)$ if and only if it is a hyperplane and it remains to show that a shard is a hyperplane if and only if it is a 
facet hyperplane of $B$.

A facet hyperplane of $B$ is necessarily basic in any rank-two subarrangement containing it, and thus it is a shard.
Conversely, suppose for the sake of contradiction that $\Sigma=H$ is a shard of $(\A,B)$ but not a facet hyperplane of $B$, and let $\B$ be 
the set of facet hyperplanes of $B$.
By Proposition~\ref{unique} there is a unique join-irreducible $J$ associated to $\Sigma$, and by Proposition~\ref{j sigma},
$J$ is the meet of all elements of $U(\Sigma)$.
But since $\Sigma$ is an entire hyperplane $H$, the region $J$ is in fact the meet of all the regions $R$ with $H\in S(R)$.
For each facet hyperplane $H'$ there is a region $R(H')$ adjacent to $B$ whose separating set is $\set{H'}$.
Thus for every $H'\in \B$ we have $H\not\in S(R(H'))$.
The join of the regions $\set{R(H'):H'\in \B}$ is $-B$.
By the antipodal anti-automorphism, we have a set of regions $\set{-R(H'):H'\in \B}$ all containing $H$ in their separating set, whose meet
is $B$, so $J=B$.
This contradicts the fact that $H\in S(J)$ and therefore no such shard $\Sigma$ exists.

We now show that (ii) is equivalent to (iii).
If $\Sigma$ is a facet hyperplane of $B$ then in particular it has no facets, so (iii) holds.
Conversely, suppose there is no facet of $\Sigma$ intersecting $J:=J_\Sigma$ in dimension $d-2$.
Let $H:=H_\Sigma$ and let $H'$ be any facet hyperplane of $J_*$.
Because $J>J_*$ we have $H\not\in S(J_*)$.
If $H\neq H'$ then since $\A$ is simplicial, the intersection $H\cap H'$ defines a $(d-2)$-dimensional face of $J$.
Let $R$ be the region with $H\in S(R)$ such that $R\cap H$ has dimension $d-1$ and $R\cap H\cap H'\cap J$ has dimension $d-2$.
Since there is no facet of $\Sigma$ intersecting $J$ in dimension $d-2$, the intersection $R\cap H$ is in $\Sigma$, so $R\in U(\Sigma)$.
Since $J$ is the unique minimal element of $U(\Sigma)$ we have $J<R$, so that $H'\not\in S(J)$ and therefore $H'\not\in S(J_*)$.
We have shown that $J_*$ is not separated from $B$ by any of its facet hyperplanes.
Therefore $J_*=B$ and $\Sigma$ is a facet hyperplane of $B$. 
\end{proof}

One can determine the poset $\Irr(\Con(\Po(\A,B)))$ computationally when $\A$ is simplicial and $\Sh(\A,B)$ is acyclic, using 
Theorem~\ref{shard} and Lemma~\ref{Cg}.
If $\Sigma_1\to\Sigma_2$, let $\A'$ be the rank-two subarrangement in which $H_1\to H_2$.
Let $D$ be some $d$-dimensional ball centered somewhere in $\Sigma_1\cap\Sigma_2$ such that $D$ intersects no hyperplanes except those 
in $\A'$.
The set $\R$ of regions intersecting $D$ forms an interval in $\Po(\A,B)$ isomorphic to $\Po(\A',B')$.
The arrow $\Sigma_1\to\Sigma_2$ corresponds to some arrow $\Sigma'_1\to\Sigma'_2$ in $\Sh(\A',B')$.
Any interval $\R$ in $\Po(\A,B)$ isomorphic to a rank-two poset of regions arises in this way.
Thus to determine $\Sh(\A,B)$, one can consider all such intervals $\R$, and given an arrow $\Sigma_1'\to\Sigma_2'$, use Lemma~\ref{Cg}
to determine the join-irreducibles $J_1$ and $J_2$ corresponding to $\Sigma_1$ and $\Sigma_2$.
For a given arrow $\Sigma_1\to\Sigma_2$, there may be more than one interval $\R$ giving rise to $\Sigma_1\to\Sigma_2$.
In Section~\ref{cox} we use this approach to give an algorithm for determining the congruence lattice of the weak order on a finite Coxeter group.

Another approach is useful when the join-irreducibles have a good combinatorial description.
One uses the combinatorial description and Lemma~\ref{polyhedra} to give an explicit description of the shards as polyhedra.
Then Lemma~\ref{arrows} provides a means of writing down the arrows in $\Sh(\A,B)$.
In Sections~\ref{Con Bn} and~\ref{Con Sn} we use this approach to determine the shard arrows for Coxeter arrangements of types B and A.

\section{Coxeter groups and weak order}
\label{cox}
In this section, we give background information on Coxeter groups and the weak order, and relate these concepts to hyperplane arrangements.
We describe an algorithm for computing $\Irr(\Con(W))$ for a finite Coxeter group $W$.
For more information on Coxeter groups and the weak order, see for example~\cite{orderings,Humphreys,hyperplane,hplanedim}.

A {\em Coxeter system} $(W,S)$ is a group $W$ given by generators $S$, and relations $s^2=1$ for all $s\in S$ and the {\em braid relations} 
$(s_1s_2)^{m(s_1,s_2)}=1$ for all $s_1\neq s_2\in S$.
Each $m(s,t)$ is an integer $\ge 2$ or is $\infty$, where $x^\infty=1$ by convention.
A Coxeter system is concisely represented by its {\em Coxeter graph} $G$, a graph whose vertex set is $S$, with edges for every pair $s_1,s_2\in S$
with $m(s_1,s_2)>2$.
The edges are labeled by $m(s_1,s_2)$, except that edge-labels $3$ are omitted.
Usually one refers to the {\em Coxeter group} $W$, implying the existence of some $S$ so that $(W,S)$ is a Coxeter system.
Important examples of Coxeter groups include finite reflection groups and Weyl groups.

Each element of $W$ can be written in many different ways as a word with letters in $S$.
A word $a$ for an element $w$ is called {\em reduced} if the length (number of letters) of $a$ is minimal among words representing~$w$.
The length of a reduced word for $w$ is called the {\em length} $l(w)$ of $w$.
Elements of $S$ are called {\em simple reflections} and any conjugate of an element of $S$ is called a {\em reflection}.
The set of reflections is denoted $T$.
When $W$ is finite, it has a {\em maximal element} $w_0$, which is maximal in length, and which is an involution in $W$.

A {\em right descent} of an element $w\in W$ is a generator $s\in S$ such that $l(w)>l(ws)$.
A {\em left descent} of $w$ is $s\in S$ such that $l(w)>l(sw)$.
The right weak order on a Coxeter group $W$ is generated by covering relations $w\covers ws$ for every $w\in W$ and every right descent $s$ 
of $w$.
There is also a left weak order, but throughout this paper, we use the right weak order and refer to it simply as the weak order.
We use the letter $W$ to denote the the partially ordered set consisting of $W$ equipped with the weak order.
An element of $W$ is join-irreducible if and only if it has a unique right descent.
Given $w\in W$, define the {\em (left) inversion set} of $w$ to be $I(w)=\set{t\in T: l(tw)<l(w)}$.
For any covering relation $u\covered w$ in weak order, there is a unique left reflection $t$ associated to the covering relation.
Namely $t=uw^{-1}$ is the unique reflection such that $u=tw$.
The (right) weak order is equivalent to containment of (left) inversion sets.
The inversion set of $w_0$ is $T$, so $w_0$ is maximal in the weak order.
The map $w\mapsto w_0ww_0$ is an automorphism of the weak order which permutes the set $S$.
The maps $w\mapsto w_0w$ and $w\mapsto ww_0$ are both anti-automorphisms.
The map $w\mapsto ww_0$ has the additional property that $I(ww_0)=T-I(w)$.

Given any finite Coxeter system $(W,S)$ and $K\subseteq S$, the subgroup generated by $K$ is called the {\em parabolic subgroup} $W_K$.
The pair $(W_K,K)$ is a finite Coxeter system, and any $w\in W$ has a unique factorization $w=w_K\cdot \!\!\phantom{.}^K\! w$ which maximizes 
$l(w_K)$ subject to the constraints that $l(w_K)+l(\!\!\phantom{.}^K\! w)=l(w)$ and that $w_K\in W_K$.
The left inversion set of $w_K$ is $I(w_K)=I(w)\cap I(w_0(K))$, where $w_0(K)$ is the maximal element of $W_K$, so that $I(w_0(K))$ is 
the set of all reflections in $W_K$.
The elements of $W_K$ are a lower interval in the weak order, so an element $\gamma\in W_K$ is join-irreducible in the weak 
order on $W_K$ if and only if it is join-irreducible in the weak order on $W$.

There is a clash of terminology inherent in studying lattice or order quotients of partial orders defined on Coxeter groups.
It is standard to take the set of elements $^K\! w$ which appear in decompositions $w=w_K\cdot \!\!\phantom{.}^K\! w$ and call this set
the {\em left quotient} $^K\! W$ of $W$ with respect to $W_K$.
There is a similar definition of right quotients $W^K$.
When $W$ is partially ordered by the weak order, $^KW$ is not a lattice quotient of $W$, but $W_K$ is.
Interestingly, when $W$ is partially ordered under the Bruhat order (which we do not define here), the situation is reversed and $^KW$ is an 
order quotient of $W$, while $W_K$ is not~\cite{Order}.

Every finite Coxeter group can be realized as a {\em reflection group}---a group generated by Euclidean reflections in some $\reals^d$.
When $W$ is realized as a reflection group, the set $T$ is exactly the set of Euclidean reflections in $W$.
Each reflection fixes a {\em reflecting hyperplane}, and our results on the weak orders arise from a study of the arrangement of reflecting 
hyperplanes associated to $W$.
Such an arrangement always simplicial and is called a {\em Coxeter arrangement}.
Every region of a Coxeter arrangement is identical by symmetry, so we choose any region~$B$ to be the base region.
Once a base region has been chosen, the elements of $W$ correspond to the regions of the associated Coxeter arrangement, and the inversion 
set of an element $w$ is the separating set of the corresponding region.
Cover relations in the weak order correspond to pairs of adjacent regions, and the left reflection associated to a cover relation fixes the 
hyperplane separating the two regions. 
The map $w\mapsto w_0ww_0$ corresponds to a Euclidean symmetry of the corresponding Coxeter arrangement, fixing the base region.
The map $w\mapsto ww_0$ corresponds to the antipodal anti-automorphism of the poset of regions.
We use $W$ to refer to the pair $(\A,B)$, as for example $\Q(W)$ and $\Sh(W)$.

The weak order is known~\cite{orderings} to be a meet-semilattice in general, and a lattice when $W$ is finite.
Caspard, Le Conte de Poly-Barbut and Morvan showed that the weak order is in fact a congruence uniform lattice~\cite{boundedref}. 
One part of their proof implies in particular that the directed graph $\Q(W)$ is acyclic and thus also $\Sh(W)$ is acyclic.

We now use the method described at the end of Section~\ref{poset} to give an algorithm for determining $\Irr(\Con(W))$ for a finite Coxeter 
group $W$.
The intervals $\R$ described there are exactly the cosets of the parabolic subgroups $W_K$ where $|K|=2$.
We may as well ignore the subsets $\set{r,s}$ with $m(r,s)=2$, because these do not contribute any shard arrows.
Let $K$ vary over all sets $\set{r,s}\subseteq S$ forming an edge in the Coxeter graph.
Let $w$ vary over all elements of $W^{\set{r,s}}$ to obtain cosets of the form $wW_{\set{r,s}}$.
For every join-irreducible $\gamma$ in $W_{\set{r,s}}$ covering $\gamma_*$ in $W_{\set{r,s}}$, determine the (left) reflection $t$ corresponding to the edge
$w\gamma\covers w\gamma_*$, and test all the join-irreducibles $\gamma'$ of $W$ with the same associated reflection to find one for which 
$(a(t)\cap I(\gamma'))\subseteq I(w\gamma)$.  
By Lemma~\ref{Cg}, this $\gamma'$ is unique and we denote it by $c(w\gamma)$.
Thus we obtain arrows of the form $c(wr)\to c(w\gamma)$ and $c(ws)\to c(w\gamma)$ for every join-irreducible $\gamma$ in $W_{\set{r,s}}$ with $l(\gamma)>1$.
Every arrow of $\Sh(W)$ arises in this way for some $\set{r,s}$ and $w$, although some arrows may arise more than once.

\section{Coxeter groups of types A and B}
\label{coxab}
The finite Coxeter groups are classified, and the infinite families are $A_n$, $B_n$, $D_n$ and $I_2(m)$.
In this section we present the usual combinatorial descriptions of the Coxeter groups $A_n$ and $B_n$, along with the corresponding 
Coxeter arrangements, and give a well-known description of the join-irreducibles of the weak order.

\subsection*{Type A (The symmetric group)}
For $n\ge 1$, the Coxeter group $A_{n-1}$ is the symmetric group $S_n$ consisting of all permutations $w$ of the set $[n]$.
We write $w=w_1w_2\cdots w_n\in S_n$, where $w_i:=w(i)$.
The generating set $S$ is the set $\set{s_i:i\in[n-1]}$, where $s_i$ is the adjacent transposition $(i,i+1)$.
The reflection set $T$ is the set of all transpositions, and the inversion set $I(w)$ of a permutation $w$ is 
$\set{(w_j,w_i):1\le i<j\le n,w_i>w_j}$.
Moving up in the weak order on $S_n$ corresponds to switching adjacent entries in a permutation so as to create an additional inversion.
The anti-automorphism $w\mapsto ww_0$ corresponds to reversing the order of entries in $w$.

The symmetric group can be realized as a reflection group so that the corresponding Coxeter arrangement $\A$ consists of the hyperplanes 
whose normal vectors are $\set{e_i-e_j:1\le j<i\le n}$, where $e_i$ is the $i$th standard basis vector of $\reals^n$.
We blur the distinction between the hyperplanes and this choice of normals, so that we refer, for example, to ``the hyperplane 
$e_i-e_j$,'' and to the inner product of two hyperplanes.
Choose $B$ to be the region such that $\br{x,e_i-e_j}<0$ for every $x$ in the interior of $B$ and for every hyperplane.
Thus for each region $R$, the separating set $S(R)$ is $\set{e_a-e_b:\br{x,e_a-e_b}>0}$ for any $x$ in the interior of $R$.
We associate to each permutation $w$ the region $R$ containing the vector $(w_1,\ldots,w_n)$.
 Thus $S(R)=\set{e_a-e_b:(b,a)\in I(w)}$.

Let $J$ be a join-irreducible region of $\Po(\A,B)$, for $(\A,B)$ as in the previous paragraph.
Then $J$ corresponds to an element $\gamma$ of $S_n$, with unique right descent $s_i$ for some $i\in[n-1]$.
That is, $\gamma_i>\gamma_{i+1}$ but $\gamma_j<\gamma_{j+1}$ for every other $j\in[n-1]$.
Let $A:=\set{\gamma_{i+1},\gamma_{i+2},\ldots,\gamma_n}$.
For any nonempty subset $A\subseteq [n]$, let $A^c:=[n]-A$ and set $m=\min A$ and $M=\max A^c$.
As long as $M>m$, we can construct a join-irreducible permutation consisting of the elements of $A^c$ in ascending order followed by 
the elements of $A$ in ascending order.
Thus join-irreducible elements of the weak order on $S_n$ correspond to nonempty subsets $A$ with $M>m$, and this correspondence is 
used frequently in later sections.
The separating set of the region $J$ is
\[S(J)=\set{e_a-e_b:a\in A^c,b\in A,a>b}.\]

\subsection*{Type B}
The Coxeter group $B_n$ is the group of {\em signed permutations}.
These are permutations $w$ of $\pm[n]:=\set{i:|i|\in[n]}$ satisfying $w(i)=-w(-i)$.
We represent a signed permutation by its {\em full notation} $w_{-n}w_{-n+1}\cdots w_{-1}w_1w_2\cdots w_n$ where 
$w_i:=w(i)$.
However, $w$ is determined by its values on the set $[n]$.
For $i\in[n-1]$, let $s_i$ be the product of transpositions $(i,i+1)(-i-1,-i)$, and let $s_0$ be the transposition $(-1,1)$.
The generating set $S$ is the set $\set{s_i:i\in[0,n-1]}$.
The reflection set $T$ is the set of all symmetric transpositions $(-a,a)$ for $a\in[n]$ together with all symmetric pairs of transpositions 
$(a,b)(-b,-a)$ for $a,b\in\pm[n]$ with $a\neq -b$.
We refer to the latter type of reflection by specifying either the transposition $(a,b)$ or the transposition $(-b,-a)$.
It is important in what follows to remember that some reflections can be named in more than one way.
We write the inversion set of a signed permutation $w$ as $\set{(w_j,w_i):i,j\in\pm[n],i<j,w_i>w_j}$,
and this expression is redundant, in that it may give some of the inversions of $w$ twice.

The group $B_n$ can be realized as a reflection group so that the corresponding Coxeter arrangement consists of the hyperplanes whose 
normal vectors are $\set{e_i:i\in[n]}\cup\set{e_i\pm e_j:1\le j<i\le n}$.
However, to avoid breaking into a  large number of cases when we determine $\Sh(B_n)$, we take a different approach.
Number the unit vectors of $\reals^{2n}$ as $e_{-n},\ldots,e_{-1},e_1,\ldots,e_n$, and consider the arrangement $\A$ 
(corresponding to $S_{2n}$) whose normals are $\set{e_i-e_j:i,j\in\pm[n],i>j}$.
For the rest of the paper, when discussing $B_n$ we assume that all subscripts $i$  in $e_i$, $x_i$, etc.\ are in $\pm[n]$, so that in particular, 
$i>j$ implies that $n\ge i>j\ge -n$ and neither $i$ nor $j$ is zero.
The arrangement corresponding to $B_n$ is the set of hyperplanes in the subspace $\set{x: x_i=-x_{-i}\mbox{ for all }i}$ of 
$\reals^{2n}$ which are contained in hyperplanes of $\A$.
We name the hyperplanes of $B_n$ by the normals to corresponding hyperplanes in $\A$.
Thus a hyperplane in $B_n$ may be referred to in more than one way.
Specifically, the hyperplane $e_i-e_j$ with $j\neq -i$ can also be referred to as $e_{-j}-e_{-i}$.
Choose $B$ to be the intersection of the base region of $\A$ (as chosen above for $S_n$) with the subspace 
$\set{x: x_i=-x_{-i}\mbox{ for all }i}$.
Thus for each region $R$, the separating set $S(R)$ is $\set{e_a-e_b:\br{x,e_a-e_b}>0}$ for any $x$ in the interior of $R$.
We associate to each signed permutation $w$ the region $R$ containing $(w_{-n},\ldots,w_{-1},w_1,\ldots,w_n)$, so that 
$S(R)=\set{e_a-e_b:(b,a)\in I(w)}$.
This expression for the separating set is redundant because it names each hyperplane in $S(R)$ in all possible ways.

A {\em signed subset} $A$ of $[n]$ is a subset of $\pm[n]$ such that for every $i\in[n]$, we have $\set{-i,i}\not\subseteq A$.
Given a nonempty signed subset $A$, let $-A:=\set{a:-a\in A}$, let $\pm A:=A\cup-A$, let the superscript ``$c$'' mean complementation in 
$\pm[n]$ and set $m:=\min A$.
Notice that the expression $-A^c$ is unambiguous, since $(-A)^c=-(A^c)$.
If $|A|=n$, set $M:=-m$ and otherwise set $M:=\max (\pm A)^c$.
Let $J$ be a join-irreducible region in $\Po(\A,B)$ where $(\A,B)$ are as in the previous paragraph.
Then $J$ corresponds to an element $\gamma\in B_n$ whose unique descent is $s_i$ for some $i\in[0,n-1]$.
In other words $\gamma_i>\gamma_{i+1}$ but $\gamma_j<\gamma_{j+1}$ for every other $j\in[0,n-1]$, so that in particular 
$\gamma_j>0$ for every $j\in[i]$.
The set $A:=\set{\gamma_{i+1},\gamma_{i+2},\ldots,\gamma_n}$ is a signed subset of $[n]$ with $M>m$.
Conversely, given a nonempty signed subset $A$ of $[n]$ with $M>m$, form a join-irreducible signed permutation consisting of the elements of 
$-A$ in ascending order, followed by $(\pm A)^c$ in ascending order, then the elements of $A$ in ascending order.
Thus join-irreducible elements of the weak order on $B_n$ correspond to nonempty signed subsets $A$ with $M>m$.
The hyperplane $H$ separating $J$ from the unique region it covers is $e_M-e_m$ and the separating set of $J$ can be expressed 
\[S(J)=\set{e_a-e_b:a>b, (a,b)\in (A^c\times A)\cup(-A\times (\pm A)^c)}.\]

The group $S_n$ is the parabolic subgroup of $B_n$ generated by $S-\set{s_0}$.
It consists of all signed permutations $w$ with $w(i)\in[n]$ for every $i\in[n]$.
Alternately, we can view $B_n$ as the sublattice of $S_{2n}$ consisting of elements fixed by the automorphism $w\mapsto w_0ww_0$.

\section{Parabolic congruences}
\label{first}
In this section, we define the notion of a parabolic congruence on $\Po(\A,B)$ in the case where $\Po(\A,B)$ is a lattice.
We define the degree of a join-irreducible, and define homogeneous congruences to be congruences generated by join-irreducibles all of the 
same degree.
The parabolic congruences are homogeneous of degree 1.

Throughout this section, $\A$ is a central hyperplane arrangement and $B$ is a base region such that $\Po(\A,B)$ is a lattice.
In \cite{BEZ} it is shown that this implies that $B$ is a simplicial region.
Recall also that if $\A$ is simplicial then $\Po(\A,B)$ is a lattice for any $B$.
Let $\B$ be the set of facet hyperplanes of the base region, and for $K\subseteq \B$ let $L_K:=\bigcap_{H\in K}H$.
The intersection $L_K\cap B$ is full dimensional in $L_K$ and is a face of $B$.
Define $\A_K:=\set{H\in\A:L_K\subseteq H}$, so that in particular $\cap_{H\in\A_K}H=L_K$.
Let $\br{H}:=\B-\set{H}$.
We omit the easy proof of the following lemma, in which we interpret the empty intersection to mean $\A$.
\begin{lemma}
\label{int}
$A_K=\bigcap_{H\in(\B-K)}A_{\br{H}}$.
\end{lemma}

For any $\A$-region R, define $R_K$ to be the $\A_K$-region containing $R$.
The separating set of $R_K$ is $S(R)\cap\A_K$.
We can also think of $R_K$ as an $\A$-region, as explained in the following lemma.
\begin{lemma}
\label{subposet}
Given any $\A$-region $R$ and any $K\subseteq \B$, there is an $\A$-region $R_K$ whose separating set is $S(R)\cap\A_K$.
\end{lemma}
\begin{proof}
By induction on $|\B-K|$.  If $K=\B$, the result is trivial because $\A=\A_K$.
Suppose $S(R)\neq(S(R)\cap\A_K)$, let $p$ be a point in the interior of $R$, and let $H\in(\B-K)$.
If $S(R)\neq(S(R)\cap\A_{\br{H}})$, let $L$ be the line $\Span(\A)\cap L_{\br{H}}$.
Since no hyperplane in $\A-\A_{\br{H}}$ contains $L$, the affine line $p+L$ intersects every hyperplane of $\A-\A_{\br{H}}$.
The intersection $L\cap B$ is a ray, and moving along $p+L$ in the direction of that ray, we eventually reach a region $R'$ separated from $B$ 
only by hyperplanes of $\A_{\br{H}}$.
Since $p+L$ intersects no hyperplane of $\A_{\br{H}}$, the separating set of this region is $S(R)\cap\A_{\br{H}}$.
If $S(R)=(S(R)\cap\A_{\br{H}})$, set $R'=R$.

In either case we have a region $R'$ whose separating set is $S(R)\cap\A_{\br{H}}$.
By induction there is a region whose separating set is $S(R')\cap\A_{K\cup\set{H}}=S(R)\cap\A_{\br{H}}\cap\A_{K\cup\set{H}}$.  
By Lemma \ref{int}, this is $S(R)\cap\A_K$.
\end{proof}

Thus it makes sense to think of $\Po(\A_K,B_K)$ as an induced subposet of $\Po(\A,B)$.
In fact it is the interval $[B,(-B)_K]$ in $\Po(\A,B)$, and furthermore it is a homomorphic image of $\Po(\A,B)$, as we now show.
\begin{prop}
\label{quotient}
$\eta_K:R\mapsto R_K$ is a lattice homomorphism.
\end{prop}
\begin{proof}
We check that $\eta_K$ is an order homomorphism.
There is a minimal element of $\Po(\A,B)$ mapping to $R_K$, namely $R_K$ itself, and the map is order preserving because 
$S(R_K)=S(R)\cap\A_K$.
The region $-B$ has the same set $\B$ of facet-hyperplanes as $B$, so we can use the same set $K$ to define a map on the anti-automorphic
poset $\Po(\A,-B)$.
The fibers of this map coincide with the fibers of $\eta$, and so in particular the fibers of $\eta$ have a unique maximal element, and
projection up to that element is order-preserving.
\end{proof}

Let $\Theta_K$ be the lattice congruence whose congruence classes are the fibers of $\eta_K$.
Recall from Section~\ref{lattice} that if $L$ is a finite lattice and $\pidown$ and $\piup$ are the projections associated to some congruence 
$\Theta$ then $\pidown(L)\cong L/\Theta$ and $\piup$ maps $\pidown(L)$ isomorphically onto $\piup(L)$.
The downward projection corresponding to the parabolic congruence $\Theta_K$ is $R\mapsto R_K$, and we denote the upward projection
by $R\mapsto R^K$.
We call the image of this upper projection $\Po(\A_K,B_K)^K$.
\begin{lemma}
\label{upper}
Let $K\subseteq \B$ and let $\Theta$ be any lattice congruence on $\Po(\A,B)$.
Then the restriction of $\Theta$ to $\Po(\A_K,B_K)$ corresponds to the restriction of $\Theta$ to $\Po(\A_K,B_K)^K$ via the isomorphism
$R\mapsto R^K$.
\end{lemma}
\begin{proof}
Suppose $R\le Q$ in $\Po(\A_K,B_K)$, so that $R^K\le Q^K$ in $\Po(\A_K,B_K)^K$.
We have $S(R^K)=S(R)\cup(\A-\A_K)$ and $S(Q^K)=S(Q)\cup(\A-\A_K)$, but $S(Q)-S(R)\subseteq\A_K$.
Thus $Q\join R^K=Q^K$ and $Q\meet R^K=R$.
If $Q\equiv R$ then $Q\join R^K\equiv R\join R^K$, so that $Q^K\equiv R^K$ and if $Q^K\equiv R^K$ then 
$Q\meet Q^K\equiv Q\meet R^K$, so that $Q\equiv R$.
\end{proof}

Our next goal is to show that when $\A$ is simplicial $\Theta_K$ has a very simple characterization as an element of $\Con(\Po(\A,B))$.
We begin with some technical lemmas.
\begin{lemma}
\label{join}
$\Theta_K$ is the join $\bigjoin_{H\in(\B-K)}\Theta_{\br{H}}$ in $\Con(\Po(\A,B))$.
\end{lemma}
\begin{proof}
We have $\Theta_K\ge\Theta_{\br{H}}$ in $\Con(\Po(\A,B))$ for every $H\in(\B-K)$ because for regions $R,Q$, if 
$S(R)\cap\A_{\br{H}}=S(Q)\cap\A_{\br{H}}$ then $S(R)\cap\A_K=S(Q)\cap\A_K$.
Let $\Phi\ge\Theta_{\br{H}}$ in $\Con(\Po(\A,B))$ for every $H\in(\B-K)$, and let $R$ be a region.
Following the argument in Lemma \ref{subposet}, we construct a sequence $R=R^1,R^2,\cdots,R^m=R_K$ such that for each $i\in[m-1]$ we 
have $R^i\equiv R^{i+1}\mod\Theta_{\br{H}}$ for some $H\in(\B-K)$.
Since $\Phi\ge\Theta_{\br{H}}$, we have $R\equiv R_K\mod\Phi$.
Therefore any other region $Q$ with $Q\equiv R\mod\Theta_K$ has $Q\equiv R\mod\Phi$, and thus $\Phi\ge\Theta_K$.
\end{proof}

\begin{lemma}
\label{parabolic lem}
Let $H_1\in(\A-\A_{\br{H}})$, let $H_2\in\A_{\br{H}}$ and let $\A'$ be the rank-two subarrangement determined by $H_1\cap H_2$.
Then $(\A'\cap\A_{\br{H}})=\set{H_2}$ and $H_2$ is basic in $\A'$.
\end{lemma}
\begin{proof}
Suppose to the contrary that $H_3\in(\A'\cap\A_{\br{H}})$ and $H_3\neq H_2$.
Then $H_1\supseteq (H_2\cap H_3)\supseteq L_{\br{H}}$, and so by definition $H_1\in\A_{\br{H}}$, contradicting the hypothesis.
Therefore $(\A'\cap\A_{\br{H}})=\set{H_2}$.
Suppose now that $H_2$ is not basic in $\A'$.
Then there exist $H_4,H_5\in(\A-\A_{\br{H}})$ with $H_4\neq H_5$ such that $H_4$ and $H_5$ are basic hyperplanes in $\A'$.
But $B\cap L_{\br{H}}=(H_2\cap B)\cap L_{\br{H}}$ has dimension $d-2$, and by Lemma \ref{basic containment}, so does
$(H_4\cap H_5\cap B)\cap L_{\br{H}}$.
In particular, $(H_4\cap H_5)\supseteq L_{\br{H}}$, so $H_4,H_5\in \A_{\br{H}}$.
This contradiction shows that $H_2$ is basic in $\A'$.
\end{proof}

\begin{lemma}
\label{shard parabolic}
The shards of $\Po(\A_{\br{H}},B_{\br{H}})$ are exactly the shards of $\Po(\A,B)$ contained in hyperplanes of $\A_{\br{H}}$.
\end{lemma}
\begin{proof}
Lemma \ref{parabolic lem} implies that no hyperplane of $\A_{\br{H}}$ is ``cut'' into shards along an intersection with a hyperplane of 
$\A-\A_{\br{H}}$.
\end{proof}

\begin{lemma}
\label{ji parabolic}
For a join-irreducible $J_1$ associated to a shard $\Sigma_1$ contained in hyperplane $H_1$, we have $H_1\in\A_{\br{H}}$ if and only if 
$S(J_1)\subseteq\A_{\br{H}}$.
\end{lemma}
\begin{proof}
If $H_1\not\in\A_{\br{H}}$, the since $H_1\in S(J_1)$ we have $S(J_1)\not\subseteq\A_{\br{H}}$.
On the other hand, if $H_1\in\A_{\br{H}}$, by Lemma \ref{shard parabolic}, we can consider $\Sigma_1$ as a shard in 
$\Po(\A_{\br{H}},B_{\br{H}})$.
So $J_1$ and $(J_1)_{\br{H}}$ are both in $U(\Sigma_1)$, with $(J_1)_{\br{H}}\le J_1$.
By Proposition~\ref{j sigma}, $J_1$ is minimal in $U(\Sigma_1)$, so $J_1=(J_1)_{\br{H}}$, or in other words $S(J_1)\subseteq\A_{\br{H}}$.
\end{proof}

Let $\B$ be the set of facet hyperplanes of $B$, and for any $H\in \B$, denote by $R(H)$ the region whose separating set is $\set{H}$.
\begin{theorem}
\label{parabolic} 
Let $\A$ be simplicial and let $K\subseteq \B$.
Then $\Theta_K$ is the unique minimal lattice congruence with $B\equiv R(H)$ for every $H\in(\B-K)$.
\end{theorem}
\begin{proof}
By Lemma~\ref{join} it is enough to prove the theorem with $K=\br{H}$ for some \mbox{$H\in \B$}. 
A region $R$ is contracted by $\eta_{\br{H}}$ if and only if $S(R)\not\subseteq\A_{\br{H}}$.
In particular, $B\equiv R(H)\mod\Theta_{\br{H}}$ and $B\not\equiv R(H')$ for $H'\in\B-\set{H}$.
We now prove by induction on $|S(J)|$ that if a join-irreducible $J$ has $S(J)\not\subseteq\A_{\br{H}}$ then $\Cg(J)\le\Cg(R(H))$ in 
$\Irr(\Con(\Po(\A,B)))$.

Suppose $J_2$ is a join-irreducible with $S(J_2)\not\subseteq\A_{\br{H}}$.
Let $\Sigma_2$ be the shard associated to $J_2$ and let $\Phi$ be the map defined after the statement of Theorem~\ref{shard}, so that
$\Phi(\Sigma_2)=\Cg(J_2)$.
Let $H_2$ be the hyperplane containing $\Sigma_2$.
Lemma~\ref{ji parabolic} says that $H_2\not\in\A_{\br{H}}$.
If $|S(J_2)|=1$ then $J_2=R(H)$.
Otherwise by Lemma~\ref{source}, $\Sigma_2$ is not a source in $\Sh(\A,B)$, so $\Sigma_2\neq H_2$.  
Let $\A'$ be a rank-two subarrangement which defines a facet-hyperplane of $\Sigma_2$ as a polyhedron in $H_2$.
By Lemma~\ref{source} we can choose $\A'$ so that $\cap\A'$ intersects $J_2$ in dimension $d-2$.
By Lemma \ref{parabolic lem}, $\A'$ contains at most one hyperplane in $\A_{\br{H}}$, so there is a basic hyperplane $H_1$ of $\A'$ in 
$\A-\A_{\br{H}}$.
Thus some shard $\Sigma_1$ in $H_1$ has $\Sigma_1\to\Sigma_2$.
Let $J_1$ be the associated join-irreducible.
Then Lemma  \ref{ji parabolic} says that $S(J_1)\not\subseteq\A_{\br{H}}$.

Since $J_2$ intersects $\cap\A'$ in dimension $d-2$ there is a region $R$ with $S(R)=(S(J_2)-\A')\cup\set{H_1}$.
But $S(J_2)$ contains $H\in\A'$ as well as some basic hyperplane of $\A'$.
Since $H$ is not basic we have $|S(J_2)\cap\A'|\ge 2$, so $|S(R)|\le |S(J_2)|-1$.
Because $R\in U(\Sigma_1)$ we have $R\ge J_{\Sigma_1}=:J_1$.
Thus $|S(J_2)|>|S(J_1)|$ so by induction $\Cg(J_1)\le\Cg(R(H))$.
Now since $\Sigma_1\to\Sigma_2$ we have $\Cg(J_2)=\Phi(\Sigma_2)\le\Phi(\Sigma_1)=\Cg(J_1)\le\Cg(R(H))$.
\end{proof}
In particular, the poset $\Irr(\Con(\Po(\A_K,B_K)))$ is isomorphic to the order filter in $\Irr(\Con(\Po(\A,B)))$ which is the complement of 
the order ideal generated by $\set{R(H):H\in(\B-K)}$.

These constructions generalize the definition of a parabolic subgroup of a finite Coxeter group.
Thus we call the congruences $\Theta_K$ {\em parabolic congruences}, and call $\A_K$ a parabolic subarrangement of $\A$.
Recall from Section~\ref{cox} the unique factorization $w=w_K\cdot \!\!\phantom{.}^K\! w$ for any $w\in W$.
Proposition~\ref{quotient} says that the map $w\mapsto w_K$ is a lattice homomorphism of weak order, a fact that also appears
in~\cite{Jedlicka}.
For a Coxeter arrangement, each fiber of this homomorphism is isomorphic to the weak order restricted to $\!\!\phantom{.}^K\! W$, the left quotient
of $W$ with respect to $W_K$, defined in Section~\ref{cox}.
However, for general simplicial arrangements the fibers need not be mutually isomorphic.
The upper projection $\piup$ of this congruence is $\piup w=w_K\cdot \!\!\phantom{.}^K\! (w_0)$.
In the language of Coxeter groups, Theorem \ref{parabolic} is the following.
\begin{corollary}
\label{weak parabolic}
Let $(W,S)$ be a Coxeter system and let $K\subseteq S$.
Then the fibers of the map $w\mapsto w_K$ constitute the smallest lattice congruence of the weak order with $1\cong s$ for every $s\in(S-K)$.
\end{corollary}

We now proceed to define the degree of a join-irreducible.
\begin{prop}
\label{support}
For any region $R$ of $\A$, there is a unique minimal $K$ (called the {\em support} of $R$) such that $R=R_K$, or equivalently $S(R)\subseteq\A_K$.
\end{prop}
\begin{proof}
Suppose $R=R_{K_1}=R_{K_2}$, or in other words 
\[S(R)=\set{H\in S(R):L_{K_1}\subseteq H}=\set{H\in S(R):L_{K_2}\subseteq H}.\]
This is equal to $\set{H\in S(R):\Span(L_{K_1}\cup L_{K_2})\subseteq H}$, which is $S\left(R_{K_1\cap K_2}\right)$, so 
$R=R_{K_1\cap K_2}$.
Since $R=R_\B$ for any $R$, the support of $R$ is the intersection of all $K\subseteq \B$ with the property that $R=R_K$.
\end{proof}
Write $\supp(R)$ for the support of $R$.
For the weak order on a Coxeter group, the support of an element $x$ is the set of generators appearing in a reduced word for~$x$.
The existence of a congruence $\Theta_{\supp(J)}$ contracting every join-irreducible whose support is not contained in $\supp(J)$ means 
in particular that none of these join-irreducibles are above $J$ in $\Irr(\Con(\Po(\A,B)))$.
Thus in $\Irr(\Con(\Po(\A,B)))$, if $J_1\le J_2$ we have $\supp(J_2)\subseteq\supp(J_1)$.
We rephrase this as a lemma for easy reference.
\begin{lemma}
\label{par lemma}
If $K\subseteq \B$, then the poset $\Irr(\Con(\Po(\A_K,B_K)))$ is an order filter in $\Irr(\Con(\Po(\A,B)))$.
\end{lemma}
Define the {\em degree} of a region $R$ to be $\deg(R)=|\supp(R)|$.
If $J_1\le J_2$ in $\Irr(\Con(\Po(\A,B)))$, we have $\deg(J_2)\le \deg(J_1)$.
A congruence $\Theta$ on $\Po(\A,B)$ is {\em homogeneous of degree $k$} if it is generated by contracting join-irreducibles of degree $k$.
A lattice homomorphism is homogeneous of degree $k$ if its associated congruence is.
The lattice homomorphisms which project onto a parabolic subgroup are the homogeneous degree-one homomorphisms.
We will see in Section \ref{examples} that other important homomorphisms are homogeneous.

When $\Po(\A,B)$ is semi-distributive, as for example when $\A$ is simplicial, and in particular when $\A$ is a Coxeter arrangement,
the antipodal anti-automorphism $\alpha$ gives rise to an automorphism of $\Irr(\Con(\Po(\A,B)))$, as explained in Section~\ref{lattice}.
If $\Theta$ is a congruence on the weak order on $W$, let $\alpha(\Theta)$ be the antipodal congruence, defined by 
$\alpha(x)\equiv \alpha(y)\mod\alpha(\Theta)$ if and only if $x\equiv y\mod\Theta$.
Then $\alpha$ induces an anti-isomorphism from $W/\Theta$ to $W/(\alpha(\Theta))$.
If $\A$ is a Coxeter arrangement, then $\alpha$ is $w\mapsto ww_0$.
\begin{proposition}
\label{dual cong}
A degree-one join-irreducible $s\in W$ is contracted by $\Theta$ if and only if it is contracted by $\alpha(\Theta)$.
Let $\gamma=srs\cdots$ be a reduced word for a degree-two join-irreducible in $W$.
Then $\gamma':=\gamma_*\cdot(w_0)_{\set{r,s}}$ is a degree-two join-irreducible with reduced word of the form $rsr\cdots$ with length 
$l(\gamma')=m(r,s)-l(\gamma)+1$ and $\gamma$ is contracted by $\Theta$ if and only if $\gamma'$ is contracted by $\alpha(\Theta)$.
\end{proposition}
\begin{proof}
An element has degree one if and only if it is some $s\in S$.
It is easily checked that $(s,sw_0)$ is a subcritical pair, so that $\alpha\circ\sigma(s)=s$.
An element $\gamma$ with a reduced word of the form $srs\cdots$ is a join-irreducible if and only if $l(\gamma)\le m(s,r)-1$, and has degree
two if and only if $l(\gamma)\ge 2$.
Suppose $\gamma$ is a degree-two join-irreducible.
In $W_{\set{s,r}}$, the pair $(\gamma,\gamma_*)$ is subcritical, and so $(\gamma,\piup(\gamma_*))=(\gamma,\gamma_*(\!\!\phantom{.}^K\! w_0))$ is subcritical 
in $W$, where $K=\set{s,r}$ and $\piup$ is the upward projection associated to the parabolic congruence $\Theta_K$.
Then 
\[\alpha\circ\sigma(\gamma)=\gamma_*(\!\!\phantom{.}^K\!w_0)\cdot w_0=\gamma_*(w_0)_K^2(\!\!\phantom{.}^K\!w_0)\cdot w_0=
\gamma_*(w_0)_K(w_0)^2=\gamma_*(w_0)_K=\gamma'.\]
We have $l(\gamma')=m(r,s)-l(\gamma_*)=m(r,s)-l(\gamma)+1$, so $2\le l(\gamma)\le m(r,s)-1$.
\end{proof}

\section{Congruences on $B_n$}
\label{Con Bn}
In this section, we apply Theorem \ref{shard} and Lemma \ref{arrows} to characterize the directed graph $\Sh(B_n)$. 
The transitive closure of $\Sh(B_n)$ is $\Irr(\Con(B_n))$, the poset of join-irreducibles of the congruence lattice of weak order on $B_n$.

Since the Coxeter group $B_n$ is the set of elements of $S_{2n}$ fixed by the automorphism $w\mapsto w_0ww_0$, one might 
expect that the congruences of $B_n$ are the symmetric order ideals in $\Irr(\Con(S_{2n}))$.
However, one can check that, for example, the relationship between $\Irr(\Con(B_2))$ and $\Irr(\Con(S_4))$ is not that simple, 
so it is necessary to determine $\Irr(\Con(B_n))$ directly from Theorem \ref{shard}.
Then, since $S_n$ is a parabolic subgroup of $B_n$, we determine $\Irr(\Con(S_n))$ as an induced subposet of $\Irr(\Con(B_n))$, 
thus avoiding repetition.

We begin by determining the arrows in $\Q(B_n)$.
Recall the realization of $B_n$ in Section~\ref{coxab}, including the fact that each hyperplane may have more than one name.
The rank-two subarrangements of size $>2$, with the basic hyperplanes underlined, are
\begin{enumerate}
\item[] $\set{\ul{e_i-e_j},\ul{e_j-e_k},e_i-e_k}\mbox{ for }i>j>k\mbox{ with }j\neq -i, k\neq -i,k\neq -j$, and
\item[] $\set{e_i-e_{-i},\ul{e_i-e_j},e_i+e_{-j},\ul{e_j-e_{-j}}}\mbox{ for }i>j>0$.
\end{enumerate}
We have arrows $e_i-e_j\to e_i-e_k$ and $e_j-e_k\to e_i-e_k$ for each $i>j>k$ with $j\neq -i, k\neq -i,k\neq -j$, as well as
the four possible arrows from the set $\set{e_i-e_j,e_j-e_{-j}}$ to the set $\set{e_i-e_{-i},e_i-e_{-j}}$ for each $i>j>0$.
Arrows of the form $e_i-e_j\to e_i-e_k$ and $e_j-e_k\to e_i-e_k$ with $j=-i$ or $k=-j$ can be re-indexed to be exactly the 
arrows of the forms $e_i-e_j\to e_i-e_{-j}$ and $e_j-e_{-j}\to e_i-e_{-j}$.
So we rewrite the collection of arrows as follows.
\[\begin{array}{ccccccl}
e_i-e_j&\to&e_i-e_k&\leftarrow&e_j-e_k&&\mbox{for }i>j>k\neq -i, \mbox{ and} \\
e_i-e_j&\to&e_i-e_{-i}&\leftarrow&e_j-e_{-j}&&\mbox{for }i>j>0.
\end{array}\]

This list is complete, but since the hyperplanes can appear under more than one name, there are other ways to name these arrows.
For the purpose of determining $a(H)\cap S(J)$ for a join-irreducible $J$ with associated hyperplane $H$, it is convenient to use this less 
redundant list.
The redundancy in the expression from Section~\ref{coxab} for $S(J)$ allows us to use the less redundant list without erroneously 
leaving hyperplanes out of $a(H)\cap S(J)$.
However, for the purpose of determining the arrows in $\Sh(B_n)$, we need to write each arrow in every possible way.
We begin with the arrows $e_i-e_j\to e_i-e_k$ and $e_j-e_k\to e_i-e_k$ for every $i>j>k\neq -i$.
Since each of the hyperplanes involved can be named in up to two different ways, both of these arrows can be named in up to four ways, for
a total of eight.
However the substitution $i\mapsto -k$, $j\mapsto -j$, $k\mapsto -i$ shows that four of the eight are redundant. 
The arrow $e_i-e_j\to e_i-e_{-i}$ can also be named $e_{-j}-e_{-i}\to e_i-e_{-i}$.
Thus these are all possible ways of naming an arrow in $\Q(B_n)$.
\[\begin{array}{cccccl}
&&e_i-e_j\\
&&\downarrow\\
e_{-j}-e_{-i}&\to&e_i-e_k&\leftarrow&e_{-k}-e_{-j}&\mbox{for }i>j>k\neq -i, \mbox{ and} \\
&&\uparrow\\
&&e_j-e_k\\
\\
\\
e_i-e_j&\to&e_i-e_{-i}&\leftarrow&e_j-e_{-j}&\mbox{for }i>j>0.\\
&&\uparrow\\
&&e_{-j}-e_{-i}\\
\end{array}\]

Let $J$ be a join-irreducible region of $B_n$ with associated subset $A$ and hyperplane $e_M-e_m$.
Recall that $S(J)=\set{e_a-e_b:a>b, (a,b)\in (A^c\times A)\cup(-A\times (\pm A)^c)}$.
The sets $a(H)\cap S(J)$ have two different forms, corresponding to two classes of signed subsets.
Class~1 is the signed subsets $A$ with $|A|=n$, or equivalently, $M=-m$.
In this case $a(H)=\set{e_M-e_b:0<b<M}\cup\set{e_a-e_{-a}:0<a<M}$.
We have $M\in -A=A^c$, and $\pm A=\pm[n]$, so that
\[a(H)\cap S(J)=\set{e_M-e_b:b\in A\cap(0,M)}\cup\set{e_a-e_{-a}:a\in A^c\cap(0,M)}.\]
Class~2 is the signed subsets with $|A|<n$, or equivalently $M\neq -m$.
In this case $a(H)=\set{e_M-e_b:m<b<M}\cup\set{e_a-e_m:m<a<M},$
and therefore
\[a(H)\cap S(J)=\set{e_M-e_b:b\in A\cap(m,M)}\cup\set{e_a-e_m:a\in A^c\cap(m,M)}.\]

The shard associated to a signed subset $A$ is the polyhedron defined by the conditions $x_i=-x_{-i}$ for all $i\in[n]$ and $x_M=x_m$, as well
as the requirement that $H'(x)\ge 0$ for each hyperplane $H'\in a(H)\cap S(J)$.
For $A$ in Class~1, since $-A=A^c$ and $m=-M$, we rewrite these conditions as:
\[\begin{array}{rccl}
x_i&=&-x_{-i}&\mbox{for all }i\in[n]\\
x_M&=&x_m\\
x_M&\ge&x_a&\mbox{for all }a\in A\cap(m,M)\\
x_a&\ge& x_M&\mbox{for all }a\in A^c\cap(m,M)
\end{array}\]
The above description also applies for $A$ in Class~2.
In the case where $A$ is in Class~1, $x_M=x_m$ implies $x_M=0$, and the fourth line of the description duplicates the third line.
Thus for $A$ in Class~1, each of these inequalities defines a facet of the shard.
For Class~2, we have the following lemmas.
In the proof of each lemma, we show that an inequality defines a facet by exhibiting a vector $x$ which does not satisfy the inequality, 
but which satisfies every other condition defining the shard.
In each case, $x$ satisfies the conditions $x_M=x_m$ and $x_i=-x_{-i}$ by construction, so these conditions are not mentioned in the 
proofs. 
\begin{lemma}
\label{Bfacet1}
Let $A$ be in Class~2, and let $a\in A\cap(m,M)$.
Then the inequality $x_M\ge x_a$ defines a facet of the associated shard $\Sigma$.
\end{lemma}
\begin{proof}
Let $x$ have coordinates:
\[x_i=\left\lbrace\begin{array}{rl}
1&\mbox{if }i=-M\mbox{ or }-m, \mbox{ or if }i\in(-A-\set{-a})\cap(-M,-m),\\
-1&\mbox{if }i=m\mbox{ or }M, \mbox{ or if }i\in(A-\set{a})\cap(m,M),\\
0&\mbox{otherwise.}
\end{array}\right.\]
This definition is not self-contradictory because $A\cap\set{-M,-m}=\emptyset$.
Because $-M,-m\not\in A$ and $m,M\not\in(m,M)$, we have $a\not\in\set{\pm m,\pm M}$, so $x_a=0>-1=x_M$.
However, if $b\in A\cap(m,M)$ and $b\neq a$, we have $x_b=-1\le-1=x_M$.
For any $b$, and in particular for $b\in A^c\cap(m,M)$ we have $x_b\ge -1=x_M$.
\end{proof}

\begin{lemma}
\label{Bfacet2}
Let $A$ be in Class~2, and let $a\in A^c\cap(m,M)$.
Then the inequality $x_a\ge x_M$ defines a facet of the associated shard $\Sigma$ if and only if one of the following holds:
\begin{enumerate}
\item[(i) ] $a\not\in\set{-M,-m}$ and $-a\not\in A\cap(m,M)$, or
\item[(ii) ]$a\in\set{-M,-m}$ and $(\pm A)^c\cap(m,M)\cap(-M,-m)=\emptyset$.
\end{enumerate}
\end{lemma}
\begin{proof}
If $a\not\in\set{-M,-m}$ and $-a\not\in A\cap(m,M)$, let $x_a=-1$, $x_{-a}=1$ and $x_i=0$ for $i\neq{\pm a}$.
By hypothesis $a\in(m,M)$, so we have $a\not\in\set{\pm m,\pm M}$ and therefore $x_M=x_m=0$.
Thus $-1=x_a<x_M$.
However, for $b\in A\cap(m,M)$ we have $b\neq \pm a$, so $x_M=0\ge 0=x_b$.
For $b\in (A^c-\set{a})\cap(m,M)$ we have $x_b=0$ or $1$, so $x_b\ge 0=x_M$.

If $a\in\set{-M,-m}$ and $(\pm A)^c\cap(m,M)\cap(-M,-m)=\emptyset$, let $x$ have coordinates:
\[x_i=\left\lbrace\begin{array}{rl}
1&\mbox{if }i=m,\mbox{ or }M,\mbox{ or if }i\in(A^c-\set{a})\cap(m,M),\\
-1&\mbox{if }i=-M,\mbox{ or }-m,\mbox{ or if }i\in(-A^c-\set{-a})\cap(-M,-m),\\
0&\mbox{otherwise.}
\end{array}\right.\]
This definition is not self-contradictory because $(m,M)\cap\set{-M,-m}=\set{a}$.
We have $x_a=-1<1=x_M$, but for any $b$, and in particular $b\in A\cap(m,M)$ we have $x_b\le 1=x_M$.
For $b\in(A^c-\set{a})\cap(m,M)$, we have $x_b=1\ge 1=x_M$.

Suppose conversely that neither condition holds.
If $a\not\in\set{-M,-m}$, this means that $-a\in A\cap(m,M)$.
Then in particular, we have $(-M,-m)\cap(m,M)\neq\emptyset$.
Since $m\neq-M$, we easily see that in fact $\set{-M,-m}\cap(m,M)\neq\emptyset$.
Furthermore, both $-M$ and $-m$ are in $A^c$, and thus we have either $x_{-M}\ge x_M$ or $x_{-m}\ge x_M$, which are both equivalent 
to $x_{-M}\ge x_M$ since $x_{-m}=x_{-M}$.
Since $-a\in A\cap(m,M)$, we also have the inequality $x_M\ge x_{-a}$.
These two inequalities imply $x_{-M}\ge x_{-a}$ or equivalently $x_a\ge x_M$.

If $a\in\set{-M,-m}$, we must have some $b\in(\pm A)^c\cap(m,M)\cap(-M,-m)$.
Since $a\in\set{-M,-m}$ and $b\in(-M,-m)$, we have $b\neq\pm a$.
We have the inequalities $x_b\ge x_M$ and $x_{-b}\ge x_M$ (or equivalently $-x_b\ge x_M$), which together imply $x_{-M}\ge x_M$, 
or in other words $x_a\ge x_M$.
Thus in either case, the inequality $x_a\ge x_M$ is not facet-defining.
\end{proof}

\begin{figure}[ht]
\caption{The shards of $B_3$.}
\label{Bpic}
\centerline{\epsfbox{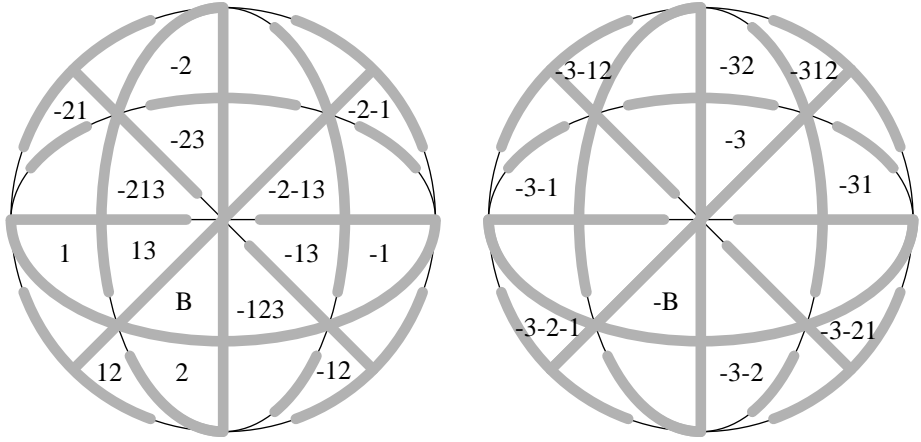}}   
\end{figure}

Figure \ref{Bpic} shows the shards of the Coxeter arrangement associated to $B_3$, with the join-irreducible regions labeled by the 
corresponding signed subset of $\set{1,2,3}$.
The signed subsets are given with set braces and commas deleted, so that for example the string $-3$$-1$2 should be interpreted as 
$\set{-3,-1,2}$.
The arrangement is represented as an arrangement of great circles on a 2-sphere.
The left drawing shows the ``northern hemisphere'' of the sphere as seen from the North Star, and the right drawing is a 180-degree rotation 
of what would be seen if the northern hemisphere were removed.
The advantage of the 180-degree rotation is that the two drawings are identical except for the labeling of the regions, and the antipodal map corresponds to translating one drawing on top of the other.
The shards are the thick gray lines, and some shards extend to both hemispheres.
The vector $e_1$ points to the right, $e_2$ points towards the top of the page, and $e_3$ points South (down into the page).

We now give a combinatorial description of the arrows in $\Sh(B_n)$ in terms of signed subsets.
Let $J_1$ and $J_2$ be join-irreducibles corresponding to signed sets $A_1$ and $A_2$ with $M_1$, $M_2$, $m_1$ and $m_2$ as defined 
in Section~\ref{coxab}.
Let $J_1$ and $J_2$ correspond to shards $\Sigma_1$ and $\Sigma_2$ contained in hyperplanes $H_1$ and $H_2$.
To simplify the statement and proof Theorem~\ref{B shard}, we first list some relevant conditions on $A_1$ and $A_2$.
We begin with conditions (q1) through (q6), in which the ``q'' indicates that these conditions determine whether $H_1\to H_2$ in $\Q(B_n)$.

\begin{tabular}{ll}
(q1)&$-m_1=M_1<M_2=-m_2$.\\
(q2)&$-m_2=M_2=M_1>m_1>0$.\\
(q3)&$M_2=M_1>m_1>m_2\neq -M_2$.\\
(q4)&$M_2>M_1>m_1=m_2\neq -M_2$.\\
(q5)&$-m_2=M_1>m_1>-M_2\neq m_2$.\\
(q6)&$-m_2>M_1>m_1=-M_2\neq m_2$.
\end{tabular}

Next we have condition (f) which is a combination of three conditions, and which depends on a parameter in $\pm[n]$.
Here the ``f'' indicates that this condition determines whether $H_1$ is a facet-defining hyperplane of $\Sigma_2$.
For $a\in\pm[n]$, say $A_2$ satisfies condition (f\,:\,$a$) if one of the following holds:

\begin{tabular}{ll}
(f1\,:\,$a$)&$a\in A_2$.\\
(f2\,:\,$a$)&$a\in A_2^c-\set{-M_2,-m_2}$ and $-a\not\in A_2\cap(m_2,M_2)$.\\
(f3\,:\,$a$)&$a\in\set{-M_2,-m_2}$ and $(\pm A_2)^c\cap(m_2,M_2)\cap(-M_2,-m_2)=\emptyset$.
\end{tabular}

\noindent
Notice that the conditions $a\in A_2$, $a\in A_2^c-\set{-M_2,-m_2}$ and $a\in\set{-M_2,-m_2}$ are mutually exclusive.

Finally, we have conditions (r1) and (r2), in which the ``r'' indicates that these conditions indicate whether $(\relint\Sigma_1)\cap\Sigma_2$ is 
non-empty.

\begin{tabular}{ll}
(r1)&$A_1\cap(m_1,M_1)=A_2\cap(m_1,M_1)$.\\
(r2)&$A_1\cap(m_1,M_1)=-A_2^c\cap(m_1,M_1)$.
\end{tabular}

\begin{theorem}
\label{B shard}
$\Sigma_1\to\Sigma_2$ if and only if one of the following combinations of conditions holds:
\begin{enumerate}
\item[1. ] {\rm(q1)} and {\rm(r1)}.
\item[2. ] {\rm(q2)} and {\rm(r1)}.
\item[3. ] {\rm(q3)}, {\rm(f\,:\,$m_1$)} and {\rm(r1)}.
\item[4. ] {\rm(q4)}, {\rm(f\,:\,$M_1$)} and {\rm(r1)}.
\item[5. ] {\rm(q5)}, {\rm(f\,:\,$\,-m_1$)} and {\rm(r2)}.
\item[6. ] {\rm(q6)}, {\rm(f\,:\,$\,-M_1$)} and {\rm(r2)}.
\end{enumerate}
\end{theorem}
\begin{proof}
We show that the theorem gives necessary and sufficient conditions on $A_1$ and $A_2$ for $\Sigma_1$ and $\Sigma_2$ to satisfy conditions 
(i), (ii) and (iii) of Lemma~\ref{arrows}.
To check condition (i), we need to use the list of arrows in $\Q(B_n)$ which includes all possible ways of naming an arrow.
We now check conditions (i) and (ii), breaking into several cases depending on which class $A_1$ and $A_2$ belong to.
\begin{enumerate}\item[]\begin{enumerate}
\case{1}{$A_1$ and $A_2$ are in Class~1.}
In this case $M_1=-m_1$ and $M_2=-m_2$.
Then $H_1\to H_2$ if and only if (q1) holds.
In this case, since $A_2$ is in Class~1, $H_1$ is a facet-defining hyperplane of $\Sigma_2$.
\case{2}{$A_1$ is in Class~2 and $A_2$ is in Class~1.}
In this case $M_1\neq-m_1$ and $M_2=-m_2$.
There is an arrow $e_{-j}-e_{-i}\to e_i-e_{-i}$ for $i>j>0$, but $M_1>0$ so $M_1\neq -j<0$.
Thus the only possibility is the arrow $e_i-e_j\to e_i-e_{-i}$ for $i>j>0$,
so $H_1\to H_2$ if and only if (q2) holds.
As in Case 1, $H_1$ is a facet-defining hyperplane of $\Sigma_2$.
\case{3}{$A_2$ is in Class~2.}
In this case $M_2\neq-m_2$, and there are four ways to have $H_1\to H_2$, corresponding to the conditions (q3) through (q6).
If (q3) holds, then by Lemmas~\ref{Bfacet1} and~\ref{Bfacet2}, $H_1$ is a facet-defining hyperplane of $\Sigma_2$ if and only if 
one of conditions (f1\,:\,$m_1$), (f2\,:\,$m_1$) or (f3\,:\,$m_1$) holds.
The conditions in Lemmas~\ref{Bfacet1} and~\ref{Bfacet2} requiring $m_1\in (m_2,M_2)$ are satisfied because (q3) holds.
Similarly, if (q4) holds, then $H_1$ is a facet-defining hyperplane of $\Sigma_2$ if and only if one of conditions (f1\,:\,$M_1$), (f2\,:\,$M_1$) or 
(f3\,:\,$M_1$) holds.
If (q5) holds, then $H_1$ is the hyperplane $e_{M_1}-e_{m_1}$, but since $-M_1=m_2$, we rename $H_1$ as $e_{-m_1}-e_{-M_1}$.
Thus in the description of $\Sigma_2$ arising from Lemma~\ref{polyhedra}, the hyperplane $H_1$ contributes some comparison between 
$x_{M_2}$ and $x_{-m_1}$.
Applying Lemmas~\ref{Bfacet1} and~\ref{Bfacet2}, we see that $H_1$ is a facet-defining hyperplane of $\Sigma_2$ if and only if 
one of conditions (f1\,:\,$-m_1$), (f2\,:\,$-m_1$) or (f3\,:\,$-m_1$) holds.
Similarly, if (q6) holds, then $H_1$ is a facet-defining hyperplane of $\Sigma_2$ if and only if one of conditions (f1\,:\,$-M_1$), (f2\,:\,$-M_1$) 
or (f3\,:\,$-M_1$) holds.
\end{enumerate}\end{enumerate}

When (q1), (q2), (q3) or (q4) holds, we describe $\Sigma_1\cap\Sigma_2$ by the following conditions: 
\[\begin{array}{rccccccl}
&&&&x_i&=&-x_{-i}&\mbox{for all }i\in[n]\\
x_{M_1}&=&x_{m_1}&=&x_{m_2}&=&x_{M_2}\\
&&&&x_{M_2}&\ge&x_a&\mbox{for all }a\in A_1\cap(m_1,M_1)\\
&&&&x_a&\ge& x_{M_2}&\mbox{for all }a\in A_1^c\cap(m_1,M_1)\\
&&&&x_{M_2}&\ge&x_a&\mbox{for all }a\in A_2\cap(m_2,M_2)\\
&&&&x_a&\ge& x_{M_2}&\mbox{for all }a\in A_2^c\cap(m_2,M_2),
\end{array}\]
which may imply that $x_{M_2}=0$.
So if $(\relint\Sigma_1)\cap\Sigma_2$ is nonempty then
\begin{eqnarray*}
A_1\cap(m_1,M_1)\cap A_2^c\cap(m_2,M_2)&=&\emptyset, \mbox{ and}\\
A_1^c\cap(m_1,M_1)\cap A_2\cap(m_2,M_2)&=&\emptyset.
\end{eqnarray*}
By (q1), (q2), (q3) or (q4), we have $(m_1,M_1)\subset(m_2,M_2)$, 
so this is equivalent to $A_1\cap(m_1,M_1)=A_2\cap(m_1,M_1)$, which is condition (r1).
Conversely, if (r1) holds, then $\Sigma_1\cap\Sigma_2=H_1\cap\Sigma_2$.
When $H_1$ is a facet-defining hyperplane of $\Sigma_2$, in particular $\Sigma_1\cap\Sigma_2$ has dimension $n-2$.
Since $H_2\not\to H_1$, $H_2$ is not a facet-defining hyperplane of $\Sigma_1$, so $\Sigma_1\cap\Sigma_2$ contains a point in 
$(\relint\Sigma_1)$.

When (q5) or (q6) holds, we can write $\Sigma_1\cap\Sigma_2$ as 
\[\begin{array}{rccccccl}
&&&&x_i&=&-x_{-i}&\mbox{for all }i\in[n]\\
x_{M_1}&=&x_{m_1}&=&-x_{m_2}&=&-x_{M_2}\\
&&&&x_a&\ge&x_{M_2}&\mbox{for all }a\in -A_1\cap(-M_1,-m_1)\\
&&&&x_{M_2}&\ge& x_a&\mbox{for all }a\in -A_1^c\cap(-M_1,-m_1)\\
&&&&x_{M_2}&\ge&x_a&\mbox{for all }a\in A_2\cap(m_2,M_2)\\
&&&&x_a&\ge& x_{M_2}&\mbox{for all }a\in A_2^c\cap(m_2,M_2),
\end{array}\]
Thus if $(\relint\Sigma_1)\cap\Sigma_2$ is nonempty then
\begin{eqnarray*}
-A_1\cap(-M_1,-m_1)\cap A_2\cap(m_2,M_2)&=&\emptyset, \mbox{ and}\\
-A_1^c\cap(-M_1,-m_1)\cap A_2^c\cap(m_2,M_2)&=&\emptyset.
\end{eqnarray*}
By (q5) or (q6) we have $(-M_1,-m_1)\subset(m_2,M_2)$, so this is equivalent to $-A_1\cap(-M_1,-m_1)=A_2^c\cap(-M_1,-m_1)$, which is
equivalent to (r2).
We finish the argument as in the case of (q1) through (q4).
\end{proof}

The poset $\Irr(\Con(B_3))$ is shown in Figure \ref{IrrConB3}.
The elements are signed subsets, with set braces and commas deleted as in Figure~\ref{Bpic}, and the antipodal 
symmetry of $\Irr(\Con(B_3))$ corresponds to reflecting the diagram through a vertical line.

\begin{figure}[ht]
\caption{$\Irr(\Con(B_3))$}
\label{IrrConB3}
\centerline{\epsfbox{IrrConB3.ps}}   
\end{figure}

\section{Congruences on $S_n$}
\label{Con Sn}
In this section, we use the the results of Section~\ref{Con Bn} and the fact that $S_n$ is a parabolic subgroup of $B_n$ to 
obtain $\Irr(\Con(S_n))$, the poset of join-irreducibles of the congruence lattice of weak order on $S_n$.
Alternately $\Irr(\Con(S_n))$ can be determined directly using a much simpler version of the method of Section \ref{Con Bn}.

Since $S_n$ is a parabolic subgroup of $B_n$, the poset $\Irr(\Con(S_n))$ is an induced subposet (in fact an order filter) of $\Irr(\Con(B_n))$,
consisting of the join-irreducibles $A$ with $m>0$, or in other words, the signed subsets with no negative elements.
Let $J_1$ and $J_2$ be join-irreducibles of $S_n$ corresponding to sets $A_1$ and $A_2$ with $M_1,M_2,m_1,m_2$ as defined in 
Section~\ref{coxab}.
Let $J_1$ and $J_2$ correspond to shards $\Sigma_1$ and $\Sigma_2$.
Conditions (q1), (q2), (q5) and (q6) in Theorem~\ref{B shard} never hold.
Also, when (q3) holds, (f1\,:\,$m_1$) must also hold, and when (q4) holds, (f2\,:\,$M_1$) must also hold.
Thus $\Sigma_1\to\Sigma_2$ if and only if (r1) holds and either (q3) or (q4) holds.
Because $m=\min A$ and $M=\max A^c$, we can rewrite these conditions as follows.

\begin{theorem}
\label{A shard}
In $\Sh(S_n)$ we have  $\Sigma_1\to\Sigma_2$ if and only if one of the following holds:
\begin{enumerate}
\item[(i) ]$A_1\cap[1,M_1)=A_2\cap[1,M_1)$ and $M_2>M_1$, or
\item[(ii) ]$A_1\cap(m_1,n]=A_2\cap(m_1,n]$ and $m_2<m_1$.
\end{enumerate}
\end{theorem}

The reader familiar with {\em root systems} (see for example~\cite{Humphreys}) may notice that the transitive closure of $\Q(S_n)$ is 
isomorphic to the {\em root poset} of the corresponding root system.
This is not true for general Coxeter groups, for example $B_2$.

\begin{figure}[ht]
\caption{The shards of $S_4$.}
\label{Apic}
\centerline{\epsfbox{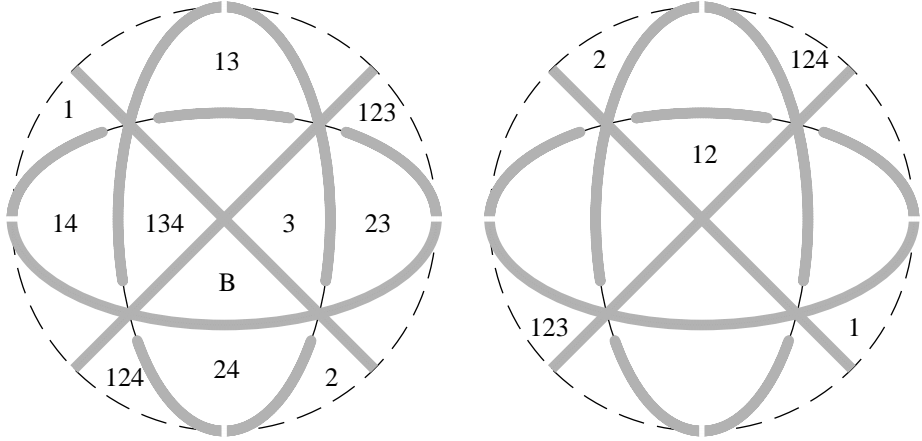}}   
\end{figure}

Figure \ref{Apic} shows the shards of the Coxeter arrangement associated to $S_4$, with the join-irreducible regions labeled by the 
corresponding subset of $\set{1,2,3,4}$.
The Coxeter arrangement corresponding to $S_n$, described in Section~\ref{coxab} is an arrangement of rank $n-1$ in $\reals^n$.
To represent $S_4$ by an arrangement in $\reals^3$, we map the normals described in Section~\ref{coxab} into $\reals^3$ by the linear map 
which fixes $e_2-e_1$ and $e_3-e_2$, sends $e_4-e_3$ to $e_2+e_1$ and sends $e_1+e_2+e_3+e_4$ to zero.
The orientation of the axes, the 180-degree rotation of the right drawing, and the other drawing conventions are the same as in 
Figure~\ref{Bpic}, except that the equatorial plane is shown as a dotted circle, indicating that it is not one of the hyperplanes in the 
arrangement.
Since the equator is not in the arrangement, some of the regions (including those labeled 1, 123, 2 and 124) intersect both hemispheres.
Also, since there are no rank-two subarrangements of size $>2$ on the equator, all of the shards which touch the equator continue through it.

Theorem~\ref{A shard} leads to a determination of the covering relations in $\Irr(\Con(S_n))$.
The proof is straightforward, and we omit the details.
\begin{theorem}
\label{A covers}
Let $A,A'\subset[n]$ represent join-irreducible elements $J$ and $J'$ of the weak order on $S_n$, let $m:=\min A$ and let $M:=\max A^c$.
Then $J$ covers $J'$ in $\Irr(\Con(S_n))$ if and only if $A'$ is one of the following:
\begin{eqnarray*}
&A-\set{M+1}&\mbox{for }M<n,\\
&(A-\set{M+1})\cup\set{M}&\mbox{for }M<n,\\
&A\cup\set{m-1}&\mbox{for }1<m \mbox{ or}\\
&(A\cup\set{m-1})-\set{m}&\mbox{for }1<m.
\end{eqnarray*} 
\end{theorem}
The degree of a join-irreducible $A$ of $S_n$ is $M-m$, and $\Irr(\Con(S_n))$ is dually ranked by $\deg(A)=M-m$, in the sense that 
$\deg(A)=\deg(A')-1$ whenever $A$ covers $A'$.
The poset $\Irr(\Con(S_4))$ is shown in Figure \ref{IrrConA3}.

\begin{figure}[ht]
\caption{$\Irr(\Con(S_4))$}
\label{IrrConA3}
\centerline{\epsfbox{IrrConA3.ps}}   
\end{figure}

\section{Examples and applications}
\label{examples}
\subsection*{The descent map}
For any simplicial arrangement $\A$ and base region $B$, let $\Delta$ be the homogeneous degree-two 
congruence on $\Po(\A,B)$ which contracts every join-irreducible of degree two or greater.
The only cover relations not contracted by $\Delta$ are those which cross facet hyperplanes of $B$.
Thus $\Po(\A,B)/\Delta$ is isomorphic to the poset of regions of the hyperplane arrangement given by the facet hyperplanes of $B$. 
Since the normals to the facet hyperplanes are a basis of $\Span(\A)$, this is a Boolean algebra.
If $\A$ is a Coxeter arrangement, the associated lattice homomorphism of the weak order is the descent map, which maps each group element 
to its left descent set.
This is because $\Delta$ contracts every covering relation whose associated left reflection is not in $S$.
Thus every element $w$ is projected down to the element $\pidown x$ which is minimal among elements whose left inversion set contains 
$I(x)\cap S$.
Since $I(x)\cap S$ is the left descent set of $x$, the element $\pidown x$ is  minimal among elements with the same descent set as $x$.
Say a lattice homomorphism is {\em of higher degree} if it contracts no join-irreducibles of degree one.
It is immediate that the descent map factors through any higher degree homomorphism.

\subsection*{The Cambrian lattices, the Tamari lattice and the cluster lattice}
In \cite{Cambrian}, the author uses a fiber-polytope construction and Theorems~\ref{A shard} and~\ref{B shard} 
to explicitly construct the {\em Cambrian lattices} of types A and B.
For any finite Coxeter group~$W$, the Cambrian lattices are a family of homogeneous degree-two lattice quotients of the weak order on $W$.
The Tamari lattices are Cambrian lattices of type A, and in type B we identify an analogous member of the family of Cambrian lattices.
Bj\"{o}rner and Wachs \cite{Nonpure II} constructed the map from $S_n$ to the Tamari lattice, and proved essentially all that was necessary to
show that it is a lattice homomorphism, but were apparently unaware of the formulation of lattice congruences in terms of order congruences.
Hugh Thomas~\cite{thomas}, working independently and approximately simultaneously, constructed the Tamari lattice of type B, and proposed a 
Tamari lattice of type D.
In types A and B (and we conjecture in all types), one of the Cambrian lattices corresponds to the {\em cluster poset}, a natural partial order 
defined in~\cite{Cambrian} on the {\em clusters}.
The clusters are the central characters in Fomin and Zelevinsky's generalized associahedra~\cite{ga}.

\subsection*{Fan Lattices}In~\cite{con_app}, it is shown that whenever $\A$ is simplicial, for any congruence $\Theta$ 
of $\Po(\A,B)$ there is a fan $\F_\Theta$ with very strong properties.
In~\cite{Cambrian}, using the results of this paper, it is shown that the fans corresponding to Cambrian congruences on Coxeter groups of 
types A and B are combinatorially isomorphic to the normal fans of associahedra.

\subsection*{Sub Hopf algebras of the Malvenuto-Reutenauer Hopf algebra}
Applying the results of the present paper, the author shows in~\cite{con_app} that certain families of congruences, consisting of a congruence on 
$S_n$ for each $n$, give rise to sub Hopf algebras of the Malvenuto-Reutenauer Hopf algebra~\cite{MR}.
In particular, we recover the well-known setup in which the Hopf algebra of quasi-symmetric functions is 
included in the Hopf algebra of planar binary trees (associated to the Tamari lattices), which is in turn included in the Malvenuto-Reutenauer 
Hopf algebra.

\subsection*{Number of congruences}
Using Theorems~\ref{A shard} and~\ref{B shard}, it is easy to write a computer program to generate the digraphs $\Sh(S_n)$ and 
$\Sh (B_n)$, up to a fairly large value of $n$.
The algorithm of Section~\ref{cox}, though less efficient, can also generate these directed graphs.
It is computationally more intense to count the order ideals in $\Irr(\Con(S_n))$ and $\Irr(\Con(B_n))$.
We used John Stembridge's {\tt posets} package \cite{posets} for Maple, as well as another program written by Stembridge, to count the 
order ideals.
The results are as follows:
\[\begin{array}{ccc|ccc}
&W&&&|\Con(W)|&\\ \hline
&S_1&&&1\\
&S_2&&&2\\
&S_3&&&7\\
&S_4&&&60\\
&S_5&&&3,\!444\\
&S_6&&&11,\!402,\!948\\
&S_7&&&129,\!187,\!106,\!461,\!769\\\hline
&B_1&&&2\\
&B_2&&&19\\
&B_3&&&8,\!368\\
&B_4&&&360,\!350,\!697,\!981\\\hline
&D_4&&&465,\!994
\end{array}\]
The sequences for $S_n$ and $B_n$ do not match any previous entries in Sloane's Encyclopedia of Integer Sequences~\cite{Sloane}.
Setting $s(n):=|\Con(S_n)|$ and $b(n) :=|\Con(B_n)|$, we can guess some approximations.
To a very good approximation for $4\le n\le 6$, we have $s(n+1)\sim s(n)^2$.  
Note that  $(3,\!444)^2 = 11,\!861,\!136$ and \[(11,\!402,\!948)^2=130,\!027,\!223,\!090,\!704.\]
The approximation $b(n+1)\sim b(n)^3$ appears to estimate the order of magnitude of $b(n)$ for these small values of $n$.
We have $19^3 = 6,\!859$ and \[(8,\!368)^3 = 585,\!956,\!012,\!032.\]

\section{Acknowledgments}
The author wishes to thank Vic Reiner and John Stembridge for helpful conversations.

\newcommand{\journalname}[1]{\textrm{#1}}
\newcommand{\booktitle}[1]{\textrm{#1}}

\end{document}